\documentclass[11pt,leqno]{article}

\usepackage{amsmath,amssymb,amsthm,amsfonts}
\usepackage{mathrsfs}
\usepackage{cite}
\usepackage{url}
\usepackage{geometry}
\usepackage{tikz}
\usetikzlibrary{patterns}

\numberwithin{equation}{section}

\newtheorem{definition}{\bf Definition}[section]
\newtheorem{theorem}[definition]{\bf Theorem}
\newtheorem{proposition}[definition]{\bf Proposition}
\newtheorem{lemma}[definition]{\bf Lemma}
\newtheorem{corollary}[definition]{\bf Corollary}
\newtheorem{remark}[definition]{\bf Remark}

\newcommand{\R}{\mathbb{R}}
\newcommand{\N}{\mathbb{N}}
\newcommand{\rmd}{\mathrm{d}}
\newcommand{\KD}{K_D}
\newcommand{\LD}{L_D}
\newcommand{\KW}{K_W}
\newcommand{\LW}{L_W}
\newcommand{\abu}{\alpha_{\rm BU}}

\allowdisplaybreaks[0]

\begin{document}

\title{Weighted pointwise estimates for the damped wave equation \\
and global solutions to Nakao's problem in three dimensions}
\author{{\large Vladimir Georgiev}\\[3mm]
    \normalsize Department of Mathematics, \\
    \normalsize University of Pisa, Largo Bruno Pontecorvo 5, 56127, Pisa, Italy, and \\
    \normalsize Faculty of Science and Engineering, Waseda University, \\
    \normalsize 3-4-1, Okubo, Shinjuku-ku, Tokyo 169-8555, Japan and \\
    \normalsize Institute of Mathematics and Informatics, Bulgarian Academy of Sciences, \\
    \normalsize Acad. Georgi Bonchev Str., Block 8, 1113 Sofia, Bulgaria \\
    \normalsize (e-mail: georgiev@dm.unipi.it)\\[4mm]
    {\large Kosuke Kita\footnote{Corresponding author}}\\[3mm]
    \normalsize Department of Mathematics, Hokkaido University, \\
    \normalsize Kita-10, Nishi-8, Kita-ku, Sapporo, Hokkaido, 060-0810, Japan\\
    \normalsize (e-mail: kosuke.kita@math.sci.hokudai.ac.jp)}
\date{}
\maketitle

\begin{abstract}
In this paper, we establish weighted pointwise estimates for the Duhamel term of the three-dimensional damped wave equation.
The estimates combine the parabolic weight \(\langle t+|x|^2\rangle\) with the wave-cone weight \(\langle t-|x|\rangle\), and thereby account for both the diffusive part of the damped wave flow and source terms with wave-like decay.
As an application, we study Nakao's problem, a weakly coupled system consisting of one semilinear damped wave equation and one semilinear wave equation.
Under explicit conditions on the exponents of nonlinear terms, we prove small-data global existence and uniqueness of mild solutions in weighted \(L^\infty\) spaces.
The result covers both the full parabolic decay and slower pointwise decay for the damped-wave component, and allows either the standard wave decay or slower pointwise decay for the wave component.
The proof combines the new Duhamel estimate with weighted estimates for the free wave equation and a contraction argument.
\end{abstract}



\section{Introduction}

In this paper, we study global existence for the following Cauchy problem for a weakly coupled system consisting of a semilinear damped wave equation and a semilinear wave equation:
\begin{equation}\label{eq:Nakao}
\left\{
\begin{aligned}
&\partial_t^2 u-\Delta u+\partial_t u=|v|^p,
&& t>0,\quad x\in\mathbb R^3,\\
&\partial_t^2 v-\Delta v=|u|^q,
&& t>0,\quad x\in\mathbb R^3,\\
&u(0,x)=u_0(x),\qquad \partial_t u(0,x)=u_1(x),\\
&v(0,x)=v_0(x),\qquad \partial_t v(0,x)=v_1(x),
\end{aligned}
\right.
\end{equation}
where \(p,q>1\).

The problem of determining the critical curve separating global existence from finite-time blow-up for \eqref{eq:Nakao} was proposed by Emeritus Professor Mitsuhiro Nakao in order to clarify how the structural difference between damped wave equations and wave equations is reflected in the nonlinear dynamics of a coupled system.
This problem is usually referred to as \emph{Nakao's problem}.

The interaction between hyperbolic propagation and parabolic dissipation has played a fundamental role in the development of the theory of nonlinear partial differential equations.
Hyperbolic equations describe finite-speed propagation and preserve information along characteristic cones, whereas parabolic equations produce smoothing and diffusion on the spatial scale \(|x|\sim t^{1/2}\).
Many important dissipative hyperbolic equations combine these two mechanisms.
A classical example is the damped wave equation.
It has finite propagation speed for every fixed time, but its large-time behavior is governed by the heat equation.
More precisely, solutions to the linear damped wave equation approach the corresponding solutions of the linear heat equation with suitable initial data.
This property is known as the \emph{diffusion phenomenon}; see, for instance, \cite{Hsiao-Liu,Nishihara1997,Chill-Haraux,Ikehata-Nishihara,Yang-Milani,Wirth}.

The emergence of parabolic behavior from hyperbolic systems with partial dissipation is not restricted to the damped wave equation.
Closely related phenomena arise in hyperbolic conservation laws with friction, relaxation systems, thermoelasticity, and other partially dissipative models.
In such systems, the dissipative variables may determine the leading large-time asymptotic profile even though the original dynamics remain hyperbolic.
At the same time, the hyperbolic part can preserve wave fronts, characteristic geometry, oscillations, or regularity-loss effects that are not visible in a purely parabolic approximation; see, for example, \cite{Hsiao-Liu,Shizuta-Kawashima,Kawashima,Wirth,Chen-Liu}.
Understanding the transfer of decay and regularity between dissipative and nondissipative components has therefore been a central issue in the past development of the subject.
It also remains important for the future study of partially dissipative systems, especially when the dissipation acts only indirectly through coupling.

The system \eqref{eq:Nakao} provides a particularly simple nonlinear model of this interaction.
Although the equation for \(u\) is hyperbolic and has finite propagation speed, its low-frequency and large-time behavior is diffusive.
By contrast, the equation for \(v\) contains no direct damping and retains the characteristic-cone structure of the free wave equation.
Moreover, the coupling acts in both directions.
The damped component is forced by the wave component through \(|v|^p\), while the undamped wave component is forced by the diffusive component through \(|u|^q\).
Consequently, wave-type decay must be transferred to the damped wave equation, and parabolic decay must simultaneously be transferred to the wave equation.

This two-way transfer distinguishes Nakao's problem from both the damped--damped and the wave--wave systems.
The mixed system cannot be regarded as a direct perturbation of either one.
Estimates adapted only to the parabolic scale do not fully exploit the localization and additional decay of the wave component away from the light cone.
Conversely, estimates adapted only to the characteristic cone do not fully capture the stronger diffusion-driven decay of the damped component.
The critical behavior of \eqref{eq:Nakao} is therefore determined by a nontrivial competition among parabolic diffusion, hyperbolic propagation, and the nonlinear feedback generated by the two exponents \(p\) and \(q\).

We recall the corresponding scalar nonlinear problems in order to make this competition more explicit.
For the semilinear heat equation
\begin{equation*}
\partial_t w-\Delta w=w^p,
\qquad t>0,\quad x\in\mathbb R^n,
\end{equation*}
the threshold is described by the Fujita index
\begin{equation*}
\alpha_{\mathrm F}(n,p)
:=
\frac{1}{p-1}-\frac n2.
\end{equation*}
The critical case corresponds to the Fujita exponent \(p=1+2/n\), and the sign of \(\alpha_{\mathrm F}(n,p)\) determines the standard blow-up/global-existence dichotomy for small nonnegative data.

In contrast, for the semilinear wave equation
\begin{equation*}
\partial_t^2 w-\Delta w=|w|^p,
\qquad t>0,\quad x\in\mathbb R^n,
\end{equation*}
the threshold is governed by the Strauss index
\begin{equation*}
\alpha_{\mathrm{St}}(n,p)
:=
\frac{p+1}{p(p-1)}-\frac{n-1}{2}.
\end{equation*}
For the classical semilinear damped wave equation with constant damping,
\begin{equation*}
\partial_t^2 w-\Delta w+\partial_t w=|w|^p,
\qquad t>0,\quad x\in\mathbb R^n,
\end{equation*}
the critical behavior is of Fujita type, reflecting the diffusion phenomenon of the damped wave flow; see, for instance, \cite{Fujita,John,Li-Zhou,Lindblad,Strauss,Todorova-Yordanov, Zhang,Ikehata-Tanizawa,Ikeda-Wakasugi,Ikeda-Ogawa, Nishihara2003,Narazaki2004,Matsumura} and the references therein.

The same contrast appears for weakly coupled systems.
For the damped--damped system
\begin{equation*}
\left\{
\begin{aligned}
&\partial_t^2 u-\Delta u+\partial_t u=|v|^p,\\
&\partial_t^2 v-\Delta v+\partial_t v=|u|^q,
\end{aligned}
\right.
\end{equation*}
the critical curve is of Fujita type and is expressed through
\begin{equation*}
\alpha_{\mathrm{DW}}(n,p,q)
:=
\max\left\{
\frac{p+1}{pq-1},
\frac{q+1}{pq-1}
\right\}
-\frac n2.
\end{equation*}
The same index appears for the corresponding weakly coupled heat system.
Indeed, when \(p=q\), one obtains
\[
\alpha_{\mathrm{DW}}(n,p,p)
=
\frac{1}{p-1}-\frac n2
=
\alpha_{\mathrm F}(n,p),
\]
which is consistent with the scalar Fujita threshold.

On the other hand, for the wave--wave system
\begin{equation*}
\left\{
\begin{aligned}
&\partial_t^2 u-\Delta u=|v|^p,\\
&\partial_t^2 v-\Delta v=|u|^q,
\end{aligned}
\right.
\end{equation*}
the critical curve is of Strauss type and is described by
\begin{equation*}
\alpha_{\mathrm W}(n,p,q)
:=
\max\left\{
\frac{p+2+q^{-1}}{pq-1},
\frac{q+2+p^{-1}}{pq-1}
\right\}
-\frac{n-1}{2}.
\end{equation*}
When \(p=q\), this reduces to
\[
\alpha_{\mathrm W}(n,p,p)
=
\frac{p+1}{p(p-1)}-\frac{n-1}{2}
=
\alpha_{\mathrm{St}}(n,p).
\]
Thus the scalar Fujita and Strauss thresholds are recovered from the corresponding coupled-system indices.
These results show that the mechanisms determining critical curves are substantially different in diffusive systems and wave systems; see \cite{Escobedo-Herrero,Narazaki,Sun-Wang,Nishihara-Wakasugi, DelSanto-Georgiev-Mitidieri,Agemi-Kurokawa-Takamura, Georgiev-Takamura-Zhou,Ikeda-Sobajima-Wakasa, Kurokawa-Takamura-Wakasa}.

Nakao's problem lies between these two well-understood classes, but it does not inherit either critical curve in a direct way.
The difficulty is not merely that the system contains two different linear propagators.
The two components have different dominant spacetime regions.
The damped component is governed at large times by the parabolic scale
\[
|x|\sim t^{1/2},
\]
whereas the undamped wave component retains significant information near the characteristic region
\[
|x|\sim t.
\]
Since the nonlinearities exchange the two components, each Duhamel term must be estimated using decay information generated by the other spacetime geometry.

This incompatibility explains why determining the exact critical curve for Nakao's problem has remained difficult, particularly in dimensions \(n\geq2\).
The available nonexistence arguments, based on test-function methods or nonlinear iteration arguments, identify regions of the \((p,q)\)-plane in which solutions with suitable positive initial data must blow up in finite time.
Such arguments, however, do not produce global solutions in the complementary region.
Conversely, standard global-existence methods for damped--damped systems generally discard part of the localization near the wave cone, while methods developed for wave--wave systems do not fully exploit the parabolic smoothing and diffusion-driven decay of the damped component.
Closing the gap between blow-up and global existence therefore requires estimates that retain both kinds of decay simultaneously.

The one-dimensional case has been settled under the corresponding assumptions, while the known results in higher dimensions have provided partial answers, successively improving the blow-up region and the corresponding upper bounds for the lifespan; see \cite{Wakasugi,Chen-Reissig,Kita-Kusaba}.
In particular, in dimension three and under suitable sign, positivity, and support assumptions on the initial data, finite-time blow-up is known when
\begin{equation*}\label{eq:known-bu-index}
\alpha_{\mathrm{BU}}
:=
\max\left\{
\frac{p+1}{pq-1}-\frac32,
\frac{q+1}{pq-1}-\frac32,
\frac{1/p+2}{pq-1}-1
\right\}
>0.
\end{equation*}
The three terms in this index reflect different decay mechanisms.
The first two have the Fujita-type form associated with the diffusive behavior of the damped component.
The final term is genuinely mixed and reflects the interaction between the characteristic propagation of the undamped wave and the diffusion-dominated behavior of the damped wave.
The presence of these different contributions already indicates that neither the damped--damped critical curve nor the wave--wave critical curve alone can describe the dynamics of \eqref{eq:Nakao}.

The present paper is devoted to the complementary direction, namely, the construction of global solutions for sufficiently small initial data.
Our purpose is to develop a weighted pointwise framework capable of resolving the parabolic scale of the damped component and the characteristic scale of the wave component at the same time.
In particular, our result establishes global existence for \(p,q>1\) satisfying
\begin{equation*}
    \max\left\{
\frac{1/p+3/2}{pq-1}-\frac12,
\frac{(3/2)p+1}{pq-1}-\frac32
\right\}
<0.
\end{equation*}
However, a gap remains between the known blow-up and global-existence regions, and the precise critical curve is still unknown.
A more detailed comparison of these regions is given after the statements of the main results, where they are summarized in the \(p\)-\(q\) plane.

The main technical difficulty is to construct a common function space that is compatible with these two different spacetime geometries.
The linear evolution of the \(u\)-component is governed at large times by the damped-wave diffusion phenomenon and has a parabolic decay scale represented by the weight
\[
\langle t+|x|^2\rangle^{-3/2}.
\]
The evolution of the \(v\)-component is governed by the free wave propagator and has a characteristic profile represented by
\[
\langle t+|x|\rangle^{-1}
\langle t-|x|\rangle^{-\beta}
\]
for a suitable \(\beta>0\).
Consequently, the source \(|v|^p\) in the damped wave equation carries decay information concentrated near the light cone, whereas the source \(|u|^q\) in the wave equation is generated on a parabolic scale.

Standard damped-wave estimates expressed only in terms of parabolic weights do not exploit the additional decay of \(|v|^p\) away from the light cone.
Conversely, standard wave estimates do not automatically retain the full diffusive decay contained in \(|u|^q\).
The central point of our argument is to incorporate a wave-cone weight into the pointwise estimate for the damped-wave Duhamel term and to combine it with the corresponding weighted estimates for the wave Duhamel term.
This permits us to control the nonlinear feedback between the two components in a common weighted \(L^\infty\) framework and to close the small-data global-existence argument.

\subsection{Notation and main results}

We introduce the notation needed to state the main results.
Throughout the paper, for \(M\in\N\) and \(\xi\in\R^M\), we set
\begin{equation*}
\langle \xi\rangle:=1+|\xi|.
\end{equation*}
We also write
\begin{equation*}
L^\infty_{t,x}:=L^\infty((0,\infty)\times\R^3).
\end{equation*}
For \(\kappa\ge 0\), put
\begin{equation*}
\|f\|_{\infty,\kappa}:=\left\|\langle\cdot\rangle^\kappa f(\cdot)\right\|_{L^\infty(\R^3)}.
\end{equation*}
Let \(D(t)\) be the solution operator associated with the damped wave equation and let \(W(t)\) be the free wave operator in three dimensions.
Their explicit formulas are recalled in Section \ref{sec:rep}.
We define
\begin{equation*}
\KD[u_0,u_1](t,x):=D(t)(u_0+u_1)(x)+\partial_tD(t)u_0(x),
\end{equation*}
\begin{equation*}
\LD[F](t,x):=\int_0^tD(t-\tau)F(\tau,x)\,\rmd\tau,
\end{equation*}
\begin{equation*}
\KW[v_0,v_1](t,x):=W(t)v_1(x)+\partial_tW(t)v_0(x),
\end{equation*}
\begin{equation*}
\LW[F](t,x):=\int_0^tW(t-\tau)F(\tau,x)\,\rmd\tau.
\end{equation*}
Then the integral formulation of \eqref{eq:Nakao} is
\begin{equation}\label{eq:Nakao-integral}
\left\{
\begin{aligned}
&u=\KD[u_0,u_1]+\LD[|v|^p],\\
&v=\KW[v_0,v_1]+\LW[|u|^q].
\end{aligned}
\right.
\end{equation}

For \(\alpha>0\), define
\begin{equation*}
Y_1^\alpha:=\left\{(\varphi,\psi)\in W^{1,\infty}(\R^3)\times L^\infty(\R^3):\|(\varphi,\psi)\|_{Y_1^\alpha}<\infty\right\},
\end{equation*}
where
\begin{equation*}
\|(\varphi,\psi)\|_{Y_1^\alpha}:=\|\varphi\|_{\infty,2\alpha}+\|\nabla\varphi\|_{\infty,2\alpha}+\|\psi\|_{\infty,2\alpha}.
\end{equation*}
For \(0<\theta\le 1\) and \(\beta\ge 0\), define
\begin{equation*}
Y_2^{\theta,\beta}:=\left\{(\varphi,\psi)\in W^{1,\infty}(\R^3)\times L^\infty(\R^3):\|(\varphi,\psi)\|_{Y_2^{\theta,\beta}}<\infty\right\},
\end{equation*}
where
\begin{equation*}
\|(\varphi,\psi)\|_{Y_2^{\theta,\beta}}:=\|\varphi\|_{\infty,\theta+\beta}+\|\nabla\varphi\|_{\infty,\theta+\beta+1}+\|\psi\|_{\infty,\theta+\beta+1}.
\end{equation*}
The corresponding solution spaces are
\begin{equation*}
X_1^\alpha:=\left\{F\in L^\infty((0,\infty)\times\R^3):\|F\|_{X_1^\alpha}<\infty\right\},
\end{equation*}
\begin{equation*}
\|F\|_{X_1^\alpha}:=\left\|\langle t+|x|^2\rangle^\alpha F(t,x)\right\|_{L^\infty_{t,x}},
\end{equation*}
and
\begin{equation*}
X_2^{\theta,\beta}:=\left\{F\in L^\infty((0,\infty)\times\R^3):\|F\|_{X_2^{\theta,\beta}}<\infty\right\},
\end{equation*}
\begin{equation*}
\|F\|_{X_2^{\theta,\beta}}:=\left\|\langle t+|x|\rangle^\theta\langle t-|x|\rangle^\beta F(t,x)\right\|_{L^\infty_{t,x}}.
\end{equation*}
A pair \((u,v)\in X_1^\alpha\times X_2^{\theta,\beta}\) is called a global mild solution of \eqref{eq:Nakao} if \((u,v)\) satisfies \eqref{eq:Nakao-integral} almost everywhere in \((0,\infty)\times\R^3\).

The first main result is a weighted Duhamel estimate for the damped wave equation.

\begin{theorem}\label{thm:damped-duhamel}
Let \(0<\alpha\le 3/2\), \(A>\alpha\), and \(\gamma>1\).
Then there exists \(C=C(\alpha,A,\gamma)>0\) such that
\begin{equation}\label{eq:damped-duhamel}
\|\LD[F]\|_{X_1^\alpha}
\le C
\left\|\langle t+|x|^2\rangle^A\langle t-|x|\rangle^\gamma F(t,x)\right\|_{L^\infty_{t,x}}.
\end{equation}
\end{theorem}

The novelty of Theorem \ref{thm:damped-duhamel} is the presence of the wave-cone weight \(\langle t-|x|\rangle^\gamma\) in a damped wave estimate.
This additional weight allows us to use the wave decay of \(v\) in the nonlinear term \(|v|^p\).
The case \(\alpha=3/2\) gives the full parabolic decay of the damped-wave component, while the range \(0<\alpha<3/2\) is needed to cover solutions with weaker damped-wave decay.

The second main result concerns global existence for Nakao's problem.

\begin{theorem}\label{thm:global-parameter}
Let \(p,q>1\), \(0<\alpha\le 3/2\), \(0<\theta\le 1\), and \(0<\beta<\theta\).
Assume
\begin{equation}\label{eq:param-cond}
\theta p>2\alpha,
\qquad
\beta p>1,
\qquad
\alpha q>1+\beta+\frac{1+\theta}{2}.
\end{equation}
If \(0<\alpha<3/2\), assume \((u_0,u_1)\in Y_1^\alpha\).
If \(\alpha=3/2\), fix \(\alpha_0>3/2\) and assume \((u_0,u_1)\in Y_1^{\alpha_0}\).
In both cases assume \((v_0,v_1)\in Y_2^{\theta,\beta}\), and define
\begin{equation*}
\mathcal N_0:=
\begin{cases}
\|(u_0,u_1)\|_{Y_1^\alpha}+\|(v_0,v_1)\|_{Y_2^{\theta,\beta}},&0<\alpha<3/2,\\
\|(u_0,u_1)\|_{Y_1^{\alpha_0}}+\|(v_0,v_1)\|_{Y_2^{\theta,\beta}},&\alpha=3/2.
\end{cases}
\end{equation*}
Then there exists \(\varepsilon_0>0\) such that if
\begin{equation*}
\mathcal N_0\le \varepsilon_0,
\end{equation*}
then \eqref{eq:Nakao} has a unique global mild solution \((u,v)\in X_1^\alpha\times X_2^{\theta,\beta}\) satisfying
\begin{equation}\label{eq:global-estimate}
\|u\|_{X_1^\alpha}+\|v\|_{X_2^{\theta,\beta}}
\le C\mathcal N_0.
\end{equation}
\end{theorem}

The first two inequalities in \eqref{eq:param-cond} ensure that the wave decay of \(v\) is sufficient to control the damped-wave Duhamel term \(\LD[|v|^p]\).
The last inequality ensures that the parabolic decay of \(u\) controls the wave Duhamel term \(\LW[|u|^q]\).
Theorem \ref{thm:global-parameter} immediately yields the following three regions.

\begin{corollary}\label{cor:strong-strong}
Let \(p>3,q>1\) and
\begin{equation}\label{eq:ss-cond}
\frac{4p-1}{pq-1}-3<0.
\end{equation}
Then there exist \(\beta\in(1/p,1)\), \(\alpha_0>3/2\), and \(\varepsilon_0>0\) such that, for any \((u_0,u_1)\in Y_1^{\alpha_0}\) and \((v_0,v_1)\in Y_2^{1,\beta}\) satisfying \(\|(u_0,u_1)\|_{Y_1^{\alpha_0}}+\|(v_0,v_1)\|_{Y_2^{1,\beta}}\le\varepsilon_0\), there exists a unique global mild solution \((u,v)\in X_1^{3/2}\times X_2^{1,\beta}\) satisfying
\begin{equation*}
\|u\|_{X_1^{3/2}}+\|v\|_{X_2^{1,\beta}}
\le C\left(\|(u_0,u_1)\|_{Y_1^{\alpha_0}}+\|(v_0,v_1)\|_{Y_2^{1,\beta}}\right).
\end{equation*}
\end{corollary}

\begin{corollary}\label{cor:weak-strong}
Let \(1<p\le 3, q>1\) and
\begin{equation}\label{eq:ws-cond}
\frac{1/p+3/2}{pq-1}-\frac12<0.
\end{equation}
Then there exist \(\alpha\in(0,3/2)\), \(\beta\in(1/p,1)\), and \(\varepsilon_0>0\) such that, for any \((u_0,u_1)\in Y_1^\alpha\) and \((v_0,v_1)\in Y_2^{1,\beta}\) satisfying \(\|(u_0,u_1)\|_{Y_1^\alpha}+\|(v_0,v_1)\|_{Y_2^{1,\beta}}\le\varepsilon_0\), there exists a unique global mild solution \((u,v)\in X_1^\alpha\times X_2^{1,\beta}\) satisfying
\begin{equation*}
\|u\|_{X_1^\alpha}+\|v\|_{X_2^{1,\beta}}
\le C\left(\|(u_0,u_1)\|_{Y_1^\alpha}+\|(v_0,v_1)\|_{Y_2^{1,\beta}}\right).
\end{equation*}
\end{corollary}

\begin{corollary}\label{cor:strong-weak}
Let \(p>3, q>1\) and
\begin{equation}\label{eq:sw-cond}
\frac{(3/2)p+1}{pq-1}-\frac32<0\le \frac{4p-1}{pq-1}-3.
\end{equation}
Then there exist \(\theta\in(0,1)\), \(\beta\in(1/p,\theta)\), \(\alpha_0>3/2\), and \(\varepsilon_0>0\) such that, for any \((u_0,u_1)\in Y_1^{\alpha_0}\) and \((v_0,v_1)\in Y_2^{\theta,\beta}\) satisfying \(\|(u_0,u_1)\|_{Y_1^{\alpha_0}}+\|(v_0,v_1)\|_{Y_2^{\theta,\beta}}\le\varepsilon_0\), there exists a unique global mild solution \((u,v)\in X_1^{3/2}\times X_2^{\theta,\beta}\) satisfying
\begin{equation*}
\|u\|_{X_1^{3/2}}+\|v\|_{X_2^{\theta,\beta}}
\le C\left(\|(u_0,u_1)\|_{Y_1^{\alpha_0}}+\|(v_0,v_1)\|_{Y_2^{\theta,\beta}}\right).
\end{equation*}
\end{corollary}

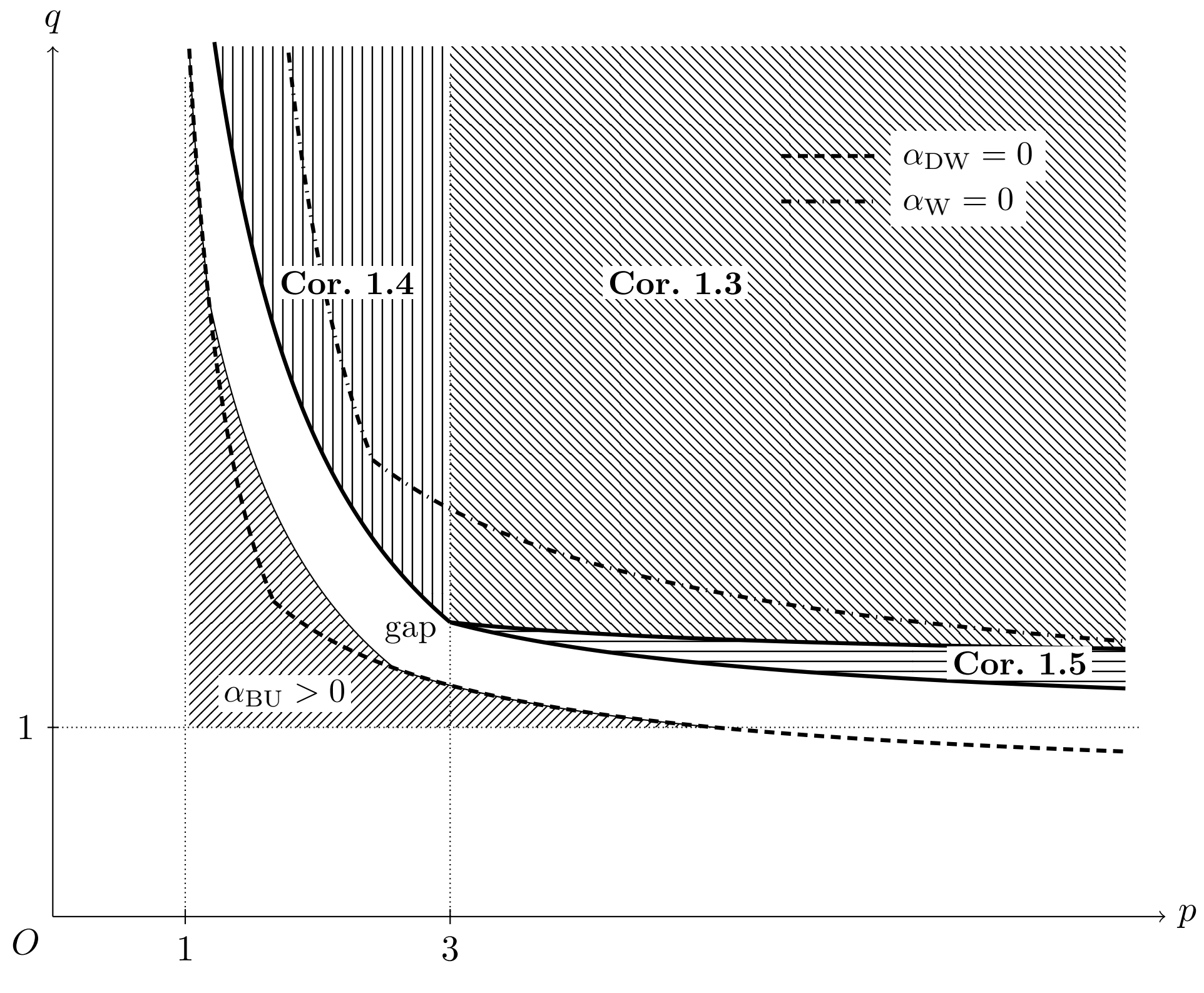
\begin{figure}[htbp]
\centering
\begin{tikzpicture}[x=1.40cm,y=2.00cm]
\let\GKsampleplot\pgfplotfunction
\def\pgfplotfunction#1#2#3{
  \edef\GKplotkey{GKpoints:#2:\detokenize{#3}}
  \ifcsname\GKplotkey\endcsname\else
    \expandafter\pgfplothandlerrecord\csname\GKplotkey\endcsname
    \GKsampleplot{#1}{#2}{#3}
    \pgfplothandlercurveto
  \fi
  \csname\GKplotkey\endcsname
}

\pgfmathsetmacro{\pzero}{(3+sqrt(41))/8}
\pgfmathsetmacro{\pone}{1+sqrt(10)/2}
\pgfmathsetmacro{\ptwo}{1+sqrt(2)}
\pgfmathsetmacro{\pBUend}{5}

\draw[->] (0,0) -- (8.4,0) node[right] {$p$};
\draw[->] (0,0) -- (0,4.6) node[above] {$q$};
\node[below left] at (0,0) {$O$};

\draw[densely dotted] (1,0) -- (1,4.45);
\draw[densely dotted] (3,0) -- (3,4.45);
\draw[densely dotted] (0,1) -- (8.2,1);

\draw (1,0.04) -- (1,-0.04) node[below] {$1$};
\draw (3,0.04) -- (3,-0.04) node[below] {$3$};
\draw (0.04,1) -- (-0.04,1) node[left] {$1$};

\path[pattern=north east lines, pattern color=black]
  (1.03,1)
  -- plot[domain=1.03:\pzero,samples=80,smooth,variable=\p]
       ({\p},{5/(3*\p-2)})
  -- plot[domain=\pzero:\pone,samples=80,smooth,variable=\p]
       ({\p},{(3*\p+1)/(\p*\p)})
  -- plot[domain=\pone:\pBUend,samples=80,smooth,variable=\p]
       ({\p},{(2*\p+5)/(3*\p)})
  -- (\pBUend,1)
  -- cycle;

\draw[thin]
  plot[domain=1.03:\pzero,samples=80,smooth,variable=\p]
       ({\p},{5/(3*\p-2)});
\draw[thin]
  plot[domain=\pzero:\pone,samples=80,smooth,variable=\p]
       ({\p},{(3*\p+1)/(\p*\p)});
\draw[thin]
  plot[domain=\pone:\pBUend,samples=80,smooth,variable=\p]
       ({\p},{(2*\p+5)/(3*\p)});

\draw[densely dashed,very thick]
  plot[domain=1.03:1.6667,samples=80,smooth,variable=\p]
       ({\p},{5/(3*\p-2)});
\draw[densely dashed,very thick]
  plot[domain=1.6667:8.1,samples=80,smooth,variable=\p]
       ({\p},{(2*\p+5)/(3*\p)});

\draw[dash dot,very thick]
  plot[domain=1.78:\ptwo,samples=80,smooth,variable=\p]
       ({\p},{(3*\p+1)/(\p*(\p-1))});
\draw[dash dot,very thick]
  plot[domain=\ptwo:8.1,samples=80,smooth,variable=\p]
       ({\p},{(\p+3+sqrt(\p*\p+10*\p+9))/(2*\p)});

\begin{scope}
\clip (1,1) rectangle (3,4.6);
\path[pattern=vertical lines, pattern color=black]
  (1.22,4.6)
  -- (3,4.6)
  -- plot[domain=3:1.22,samples=80,smooth,variable=\p]
       ({\p},{(4*\p+2)/(\p*\p)})
  -- cycle;
\end{scope}

\begin{scope}
\clip (3,1) rectangle (8.1,4.6);
\path[pattern=north west lines, pattern color=black]
  (3,4.6)
  -- (8.1,4.6)
  -- (8.1,{(4*8.1+2)/(3*8.1)})
  -- plot[domain=8.1:3,samples=80,smooth,variable=\p]
       ({\p},{(4*\p+2)/(3*\p)})
  -- cycle;
\end{scope}

\begin{scope}
\clip (3,1) rectangle (8.1,4.6);
\path[pattern=horizontal lines, pattern color=black]
  plot[domain=3:8.1,samples=200,smooth,variable=\p]
       ({\p},{(4*\p+2)/(3*\p)})
  -- plot[domain=8.1:3,samples=200,smooth,variable=\p]
       ({\p},{1+5/(3*\p)})
  -- cycle;
\end{scope}

\draw[very thick]
  plot[domain=1.22:3,samples=200,smooth,variable=\p]
       ({\p},{(4*\p+2)/(\p*\p)});

\draw[very thick]
  plot[domain=3:8.1,samples=200,smooth,variable=\p]
       ({\p},{(4*\p+2)/(3*\p)});

\draw[very thick]
  plot[domain=3:8.1,samples=200,smooth,variable=\p]
       ({\p},{1+5/(3*\p)});

\node[fill=white, inner sep=1.5pt] at (2.22,3.35)
{\bf\small Cor.\ \ref{cor:weak-strong}};
\node[fill=white, inner sep=1.5pt] at (4.7,3.35)
{\bf\small Cor.\ \ref{cor:strong-strong}};
\node[fill=white, inner sep=1.5pt] at (7.3,1.34)
{\bf\small Cor.\ \ref{cor:strong-weak}};
\node[fill=white, inner sep=1.5pt] at (1.75,1.18)
{\bf\small $\abu>0$};
\node at (2.7,1.50) {\small gap};

\begin{scope}[shift={(5.50,4.50)}]
\draw[densely dashed,very thick] (0,-0.48) -- (0.72,-0.48);
\node[fill=white,right] at (0.82,-0.48) {\small $\alpha_{\rm DW}=0$};

\draw[dash dot,very thick] (0,-0.72) -- (0.72,-0.72);
\node[fill=white,right] at (0.82,-0.72) {\small $\alpha_{\rm W}=0$};
\end{scope}

\end{tikzpicture}
\caption{Global existence and nonexistence for \eqref{eq:Nakao}}
\label{fig:pq-plane}
\end{figure}

\begin{remark}
The three global-existence regions shown in Figure~\ref{fig:pq-plane} correspond to Corollaries~\ref{cor:strong-strong}, \ref{cor:weak-strong}, and \ref{cor:strong-weak}.
Corollary~\ref{cor:weak-strong} treats the range \(1<p\leq 3\) by allowing the damped-wave component to have weaker than full parabolic decay.
For \(p>3\), Corollary~\ref{cor:strong-strong} gives global solutions with the full parabolic decay of the damped-wave component and the standard pointwise decay of the wave component.
Corollary~\ref{cor:strong-weak} extends this region downward by retaining the full parabolic decay of the damped-wave component while allowing weaker pointwise decay of the wave component.

Combining the three corollaries, we obtain global existence for \(p,q>1\) satisfying
\begin{equation*}
\max\left\{
\frac{1/p+3/2}{pq-1}-\frac{1}{2},
\frac{(3/2)p+1}{pq-1}-\frac{3}{2}
\right\}<0.
\end{equation*}
This boundary lies strictly below the critical curve \(\alpha_{\mathrm W}(3,p,q)=0\) for the corresponding wave--wave system.
In particular, our result gives global solutions in a nonempty part of the blow-up side of the wave--wave critical curve, which reflects the effect of the damping acting on the \(u\)-component.
Nevertheless, a gap remains between the known blow-up and global-existence regions for Nakao's problem, and the precise critical curve is still unknown.
\end{remark}

The paper is organized as follows.
In Section \ref{sec:rep}, we recall representation formulas and collect the auxiliary weighted \(L^\infty\) estimates for the wave and damped wave equations.
In Section \ref{sec:proof-damped}, we prove Theorem \ref{thm:damped-duhamel}.
Section \ref{sec:global} is devoted to the proof of Theorem \ref{thm:global-parameter} and Corollaries \ref{cor:strong-strong}--\ref{cor:strong-weak}.

\section{Preliminaries}\label{sec:rep}

We first recall the representation formulas used throughout the paper.
For the free wave equation in \(\R^3\), we use the propagator
\begin{equation*}
W(t)f(x):=\frac{t}{4\pi}\int_{\mathbb S^2}f(x+t\omega)\,\rmd\sigma(\omega).
\end{equation*}
This is Kirchhoff's formula.
For the damped wave equation, the corresponding propagator is denoted by \(D(t)\).
In three space dimensions, it is represented as
\begin{equation}\label{eq:D-def}
D(t)=e^{-t/2}W(t)+S(t),
\end{equation}
where
\begin{equation}\label{eq:S-def}
S(t)f(x):=\frac{e^{-t/2}}{4\pi}\int_0^t\int_{\mathbb S^2}
\frac{I_1\left(\frac{\sqrt{t^2-\rho^2}}{2}\right)}{2\sqrt{t^2-\rho^2}}\rho^2 f(x+\rho\omega)\,\rmd\sigma(\omega)\,\rmd\rho.
\end{equation}
Here \(I_\nu\) denotes the modified Bessel function of order \(\nu\),
\begin{equation*}
I_\nu(r):=\sum_{n=0}^\infty \frac{1}{n!\Gamma(n+\nu+1)}\left(\frac r2\right)^{2n+\nu},
\end{equation*}
where
\begin{equation*}
\Gamma(r):=\int_0^\infty s^{r-1}e^{-s}\,\rmd s.
\end{equation*}

\subsection{Auxiliary weighted \(L^\infty\) estimates}\label{sec:linear}

We collect the linear estimates used later.
These are not the new estimates established in this paper; the new point is the damped-wave Duhamel estimate proved in Section \ref{sec:proof-damped}.
The estimates below are included here to make the fixed point argument in Section \ref{sec:global} self-contained at the level of notation and exponents.
The constants in this subsection depend only on the displayed exponents.

\begin{proposition}[\cite{John,Georgiev-Heiming-Kubo}]\label{prop:wave-estimates}
Let \(m>2\) and \(0<\beta\le m-1\).
Then
\begin{equation}\label{eq:wave-hom-strong}
\|\KW[v_0,v_1]\|_{X_2^{1,\beta}}
\le C\|(v_0,v_1)\|_{Y_2^{1,\beta}},
\end{equation}
\begin{equation}\label{eq:wave-inhom-strong}
\|\LW[F]\|_{X_2^{1,\beta}}
\le C\left\|\langle t+|x|\rangle\langle t-|x|\rangle^\beta\langle x\rangle^mF(t,x)\right\|_{L^\infty_{t,x}}.
\end{equation}
\end{proposition}

For later use, we prove the following weak decay estimates for the wave equation.
\begin{proposition}\label{prop:wave-weak}
Let \(0<\theta<1\) and \(0<\beta<\theta\).
Then
\begin{equation}\label{eq:wave-hom-weak}
\|\KW[v_0,v_1]\|_{X_2^{\theta,\beta}}
\le C\|(v_0,v_1)\|_{Y_2^{\theta,\beta}},
\end{equation}
\begin{equation}\label{eq:wave-inhom-weak}
\|\LW[F]\|_{X_2^{\theta,\beta}}
\le C\left\|\langle t+|x|\rangle^\theta\langle t-|x|\rangle^{1-\theta+\beta}\langle x\rangle^{1+\theta}F(t,x)\right\|_{L^\infty_{t,x}}.
\end{equation}
\end{proposition}

\begin{proof}
We prove both estimates directly from the Kirchhoff formula.
For a measurable function \(h\), set
\begin{equation*}
M_t h(x):=\frac{1}{4\pi}\int_{\mathbb S^2}h(x+t\omega)\,\rmd\sigma(\omega).
\end{equation*}
Then
\begin{equation}\label{eq:Kirchhoff-weak}
\KW[v_0,v_1](t,x)=M_t v_0(x)+tM_t(\omega\cdot\nabla v_0)(x)+tM_t v_1(x).
\end{equation}
We first recall the elementary one-dimensional estimates used below.
If \(0<\beta<\theta<1\), then there exists \(C>0\) such that, for all \(t\ge0\) and \(r>0\),
\begin{equation}
\label{eq:weak-wave-1d-a}
\frac{\langle t+r\rangle^\theta\langle t-r\rangle^\beta}{r}
\int_{|t-r|}^{t+r}\langle\lambda\rangle^{-(\theta+\beta)}\,\rmd\lambda
\le C,
\end{equation}
and, for all \(t>0\) and \(r>0\),
\begin{equation}
\label{eq:weak-wave-1d-b}
\frac{\langle t+r\rangle^\theta\langle t-r\rangle^\beta}{tr}
\int_{|t-r|}^{t+r}\lambda\langle\lambda\rangle^{-(\theta+\beta)}\,\rmd\lambda
\le C.
\end{equation}

We first prove \eqref{eq:weak-wave-1d-a}.
If \(t\ge2r\), then \(|t-r|\ge t/2\), and hence \(\langle t+r\rangle\lesssim \langle\lambda\rangle\) on the interval \([|t-r|,t+r]\).
Therefore
\begin{equation*}
\int_{|t-r|}^{t+r}\langle\lambda\rangle^{-(\theta+\beta)}\,\rmd\lambda
\lesssim r\langle t+r\rangle^{-(\theta+\beta)}.
\end{equation*}
Since \(\langle t-r\rangle\le \langle t+r\rangle\), this gives \eqref{eq:weak-wave-1d-a}.
If \(t< 2r\), then \(3r>t+r\).
Assume that \(t+r>2\).
Then
\begin{equation*}
\int_{|t-r|}^{t+r}\langle\lambda\rangle^{-(\theta+\beta)}\,\rmd\lambda
\lesssim
\langle t-r \rangle^{-\beta}
\int_{0}^{t+r}\langle\lambda\rangle^{-\theta}\,\rmd\lambda
\lesssim
\langle t-r \rangle^{-\beta} \langle t+r\rangle^{1-\theta}
\lesssim r\langle t-r \rangle^{-\beta} \langle t+r\rangle^{-\theta}.
\end{equation*}
If instead \(t+r\le 2\), since \(\langle t-r\rangle\le\langle t+r\rangle\lesssim 1\), we have
\begin{equation*}
    \frac{\langle t+r\rangle^\theta\langle t-r\rangle^\beta}{r}
\int_{|t-r|}^{t+r}\langle\lambda\rangle^{-(\theta+\beta)}\,\rmd\lambda
\lesssim \frac{1}{r} \int_{|t-r|}^{t+r}\,\rmd\lambda
\lesssim 1.
\end{equation*}
This proves \eqref{eq:weak-wave-1d-a}.

Next we prove \eqref{eq:weak-wave-1d-b}.
Suppose first that \(t\ge 2r\) or \(r\ge 2t\).
Then \(3|t-r|\ge t+r\).
Hence
\begin{equation*}
\int_{|t-r|}^{t+r}
\lambda\langle\lambda\rangle^{-(\theta+\beta)}\,\rmd\lambda
\lesssim
\min\{t,r\}(t+r)\langle t+r\rangle^{-(\theta+\beta)}
\lesssim
tr\langle t+r\rangle^{-(\theta+\beta)}.
\end{equation*}
Since \(\langle t-r\rangle\le\langle t+r\rangle\), this gives \eqref{eq:weak-wave-1d-b} in this case.
It remains to consider the case \(t/2<r<2t\).
In this case, we have \(tr\gtrsim (t+r)^2\).
Since \(\theta+\beta<2\), we have
\begin{equation*}
\int_{|t-r|}^{t+r}
\lambda\langle\lambda\rangle^{-(\theta+\beta)}\,\rmd\lambda
\le
\int_0^{t+r}
\lambda\langle\lambda\rangle^{-(\theta+\beta)}\,\rmd\lambda
\lesssim
(t+r)^2\langle t+r\rangle^{-(\theta+\beta)}
\lesssim tr\langle t+r\rangle^{-(\theta+\beta)}.
\end{equation*}
Similarly, \(\langle t-r\rangle\le \langle t+r\rangle\) implies \eqref{eq:weak-wave-1d-b}.

Let \(|x|>0\).
By the spherical mean identity,
\begin{equation}\label{eq:spherical-mean-bound}
\int_{\mathbb S^2}H(|x+\rho\omega|)\,\rmd\sigma(\omega)
= \frac{2\pi}{|x|\rho}
\int_{||x|-\rho|}^{|x|+\rho}\lambda H(\lambda)\,\rmd\lambda
\end{equation}
holds for every nonnegative radial continuous function \(H\).
Applying \eqref{eq:spherical-mean-bound} with \(\rho=t\), we obtain
\begin{equation*}
|tM_t v_1(x)|+|tM_t(\omega\cdot\nabla v_0)(x)|
\lesssim
\frac{\|(v_0,v_1)\|_{Y_2^{\theta,\beta}}}{|x|}
\int_{|t-|x||}^{t+|x|}\langle\lambda\rangle^{-(\theta+\beta)}\,\rmd\lambda.
\end{equation*}
Thus \eqref{eq:weak-wave-1d-a} controls these two terms.
Similarly,
\begin{equation*}
|M_t v_0(x)|
\lesssim
\frac{\|(v_0,v_1)\|_{Y_2^{\theta,\beta}}}{t|x|}
\int_{|t-|x||}^{t+|x|}\lambda\langle\lambda\rangle^{-(\theta+\beta)}\,\rmd\lambda
\end{equation*}
for \(t>0\), and \eqref{eq:weak-wave-1d-b} applies.
The cases \(t=0\) or \(x=0\) follow directly from the definitions of \(M_t\).
Combining these estimates with \eqref{eq:Kirchhoff-weak} proves \eqref{eq:wave-hom-weak}.

Next we prove \eqref{eq:wave-inhom-weak}.
Put
\begin{equation*}
M:=\left\|\langle t+|x|\rangle^\theta\langle t-|x|\rangle^{1-\theta+\beta}\langle x\rangle^{1+\theta}F(t,x)\right\|_{L^\infty_{t,x}}.
\end{equation*}
For \(|x|>0\), the representation of \(W(t-s)\) and \eqref{eq:spherical-mean-bound} yield
\begin{align}
|\LW[F](t,x)|
&\lesssim
\frac{M}{|x|}\int_0^t\int_{||x|-(t-s)|}^{|x|+t-s}
\lambda\langle s+\lambda\rangle^{-\theta}
\langle s-\lambda\rangle^{-(1-\theta+\beta)}
\langle\lambda\rangle^{-(1+\theta)}\,\rmd\lambda\,\rmd s. \label{eq:weak-wave-duhamel-reduced}
\end{align}
It suffices to show that
\begin{align}
\frac{1}{|x|}\int_0^t\int_{||x|-(t-s)|}^{|x|+t-s}
\lambda\langle s+\lambda\rangle^{-\theta}
\langle s-\lambda\rangle^{-(1-\theta+\beta)}
\langle\lambda\rangle^{-(1+\theta)}\,\rmd\lambda\,\rmd s
\lesssim\langle t+|x|\rangle^{-\theta}\langle t-|x|\rangle^{-\beta}. \label{eq:weak-wave-duhamel-1d}
\end{align}

First, assume that \(|x|\le (t+|x|)/4\).
Then \(|x|\le t/3\), and hence \(\langle t\rangle\simeq\langle t+|x|\rangle\simeq\langle t-|x|\rangle\).
Moreover, for \(\left||x|-(t-s)\right|\le \lambda\le |x|+t-s\), we have \(t-|x|\le s+\lambda\le t+|x|\) and hence \(\langle t\rangle\simeq\langle s+\lambda\rangle\).
Therefore, for \(s\in[0,t]\), it suffices to estimate
\begin{equation*}
\frac{1}{|x|}
\int_{\left||x|-(t-s)\right|}^{|x|+t-s}
\lambda
\langle s-\lambda\rangle^{-(1-\theta+\beta)}
\langle\lambda\rangle^{-(1+\theta)}
\,\rmd\lambda
\end{equation*}

Put \(\delta=1-\theta+\beta\).
We shall use the following elementary estimate: for every \(0\le s\le t\),
\begin{equation}\label{eq:wave-weak-inner-average}
\begin{aligned}
&\frac1{|x|}
\int_{\left||x|-(t-s)\right|}^{|x|+t-s}
\lambda
\langle\lambda\rangle^{-(1+\theta)}
\langle s-\lambda\rangle^{-\delta}
\,\rmd\lambda  \\
&\qquad\le
C
\langle t-s\rangle^{-\theta}
\langle 2s-t\rangle^{-\delta}.
\end{aligned}
\end{equation}
We prove \eqref{eq:wave-weak-inner-average}.
If \(t-s\ge 2|x|\), then the integration interval is
\begin{equation*}
[t-s-|x|,\,t-s+|x|].
\end{equation*}
Thus \(\lambda\ge (t-s)/2\) on this interval, and hence
\begin{equation*}
\lambda\langle\lambda\rangle^{-(1+\theta)}
\le C\langle t-s\rangle^{-\theta}.
\end{equation*}
With \(y=s-\lambda\), we obtain
\begin{equation*}
\begin{aligned}
\frac1{|x|}
\int_{t-s-|x|}^{t-s+|x|}
\lambda
\langle\lambda\rangle^{-(1+\theta)}
\langle s-\lambda\rangle^{-\delta}
\,\rmd\lambda
\le
C\langle t-s\rangle^{-\theta}
\frac1{|x|}
\int_{2s-t-|x|}^{2s-t+|x|}
\langle y\rangle^{-\delta}
\,\rmd y.
\end{aligned}
\end{equation*}
Since \(0<\delta<1\), the elementary averaging estimate
\begin{equation*}
\frac1h\int_{A-h}^{A+h}\langle y\rangle^{-\delta}\,\rmd y
\le C\langle A\rangle^{-\delta},
\qquad A\in\mathbb R,\quad h>0,
\end{equation*}
gives \eqref{eq:wave-weak-inner-average} in this case.

If \(t-s\le |x|/2\), then the integration interval is
\begin{equation*}
[|x|-(t-s),\,|x|+t-s].
\end{equation*}
On this interval, \(\lambda\ge |x|/2\), and hence
\begin{equation*}
\lambda\langle\lambda\rangle^{-(1+\theta)}
\le C\langle |x|\rangle^{-\theta}
\le C\langle t-s\rangle^{-\theta}.
\end{equation*}
Moreover, since \(|x|\le t/3\) and \(t-s\le |x|/2\), we have
\begin{equation*}
 t-2|x|\le s-\lambda\le t
 \qquad\text{and}\qquad
 \frac{2t}{3}\le 2s-t\le t.
\end{equation*}
Therefore
\begin{equation*}
\langle s-\lambda\rangle^{-\delta}
\le C\langle 2s-t\rangle^{-\delta}
\end{equation*}
on the integration interval.
Since the length of the interval is \(2(t-s)\), we obtain
\begin{equation*}
\begin{aligned}
&\frac1{|x|}
\int_{|x|-(t-s)}^{|x|+t-s}
\lambda
\langle\lambda\rangle^{-(1+\theta)}
\langle s-\lambda\rangle^{-\delta}
\,\rmd\lambda  \\
&\qquad\le
C\frac{t-s}{|x|}
\langle t-s\rangle^{-\theta}
\langle 2s-t\rangle^{-\delta}
\le
C\langle t-s\rangle^{-\theta}
\langle 2s-t\rangle^{-\delta}.
\end{aligned}
\end{equation*}
Thus \eqref{eq:wave-weak-inner-average} follows in this case as well.

It remains to consider the case
\begin{equation*}
\frac{|x|}{2}<t-s<2|x|.
\end{equation*}
Then \(1/|x|\le C/(t-s)\), and the integration interval is contained in
\begin{equation*}
[0,\,3(t-s)].
\end{equation*}
Hence it suffices to prove
\begin{equation}\label{eq:third-inner-estimate}
\begin{aligned}
&\frac1{t-s}
\int_0^{3(t-s)}
\lambda
\langle\lambda\rangle^{-(1+\theta)}
\langle 2s-t+(t-s)-\lambda\rangle^{-\delta}
\,\rmd\lambda \\
&\qquad\le
C\langle t-s\rangle^{-\theta}
\langle 2s-t\rangle^{-\delta}.
\end{aligned}
\end{equation}
If \(|2s-t|\ge 4(t-s)\), then for \(0\le\lambda\le 3(t-s)\),
\begin{equation*}
|2s-t+(t-s)-\lambda|
\ge |2s-t|-|(t-s)-\lambda|
\ge \frac12|2s-t|.
\end{equation*}
Therefore
\begin{equation*}
\langle 2s-t+(t-s)-\lambda\rangle^{-\delta}
\le C\langle 2s-t\rangle^{-\delta}.
\end{equation*}
Since
\begin{equation*}
\frac1{t-s}
\int_0^{3(t-s)}
\lambda\langle\lambda\rangle^{-(1+\theta)}\,\rmd\lambda
\le C\langle t-s\rangle^{-\theta},
\end{equation*}
we obtain \eqref{eq:third-inner-estimate} in this subcase.

It remains to consider the subcase \(|2s-t|<4(t-s)\).
If \(t-s\le1\), then the left-hand side of \eqref{eq:third-inner-estimate} is bounded by a constant, whereas the right-hand side is bounded below by a positive constant.
Thus we may assume \(t-s>1\).
By the change of variables \(\lambda=(t-s)\eta\), the left-hand side of \eqref{eq:third-inner-estimate} is bounded by
\begin{equation*}
C(t-s)^{-\theta-\delta}
\int_0^3
\eta^{-\theta}
\left|\frac{2s-t}{t-s}+1-\eta\right|^{-\delta}
\,\rmd\eta.
\end{equation*}
Here
\begin{equation*}
\frac{2s-t}{t-s}+1>\frac12,
\qquad
\left|\frac{2s-t}{t-s}\right|<4.
\end{equation*}
The first inequality follows from \(|x|/2<t-s<2|x|\) and \(|x|\le t/3\).
Hence the last integral is bounded uniformly, since the possible singular point of the second factor stays in a compact subset of \((0,\infty)\) and \(0<\theta<1\), \(0<\delta<1\).
Consequently, the left-hand side of \eqref{eq:third-inner-estimate} is bounded by
\begin{equation*}
C(t-s)^{-\theta-\delta}.
\end{equation*}
Since \(|2s-t|<4(t-s)\), we have
\begin{equation*}
\langle t-s\rangle^{-\delta}
\le C\langle 2s-t\rangle^{-\delta},
\end{equation*}
and \eqref{eq:third-inner-estimate} follows.
This proves \eqref{eq:wave-weak-inner-average}.

Using \eqref{eq:wave-weak-inner-average} and
\begin{equation*}
\langle s+\lambda\rangle^{-\theta}
\le C\langle t\rangle^{-\theta},
\end{equation*}
we obtain
\begin{equation*}
\begin{aligned}
&\frac1{|x|}
\int_{\left||x|-(t-s)\right|}^{|x|+t-s}
\lambda
\langle s+\lambda\rangle^{-\theta}
\langle s-\lambda\rangle^{-(1-\theta+\beta)}
\langle\lambda\rangle^{-(1+\theta)}
\,\rmd\lambda  \\
&\qquad\le
C
\langle t\rangle^{-\theta}
\langle t-s\rangle^{-\theta}
\langle 2s-t\rangle^{-(1-\theta+\beta)}.
\end{aligned}
\end{equation*}

We also use the following one-dimensional convolution estimate:
\begin{equation}\label{eq:weak-wave-convolution}
\int_0^t
\langle 2s-t\rangle^{-(1-\theta+\beta)}
\langle t-s\rangle^{-\theta}
\,\rmd s
\le C\langle t\rangle^{-\beta}.
\end{equation}
The assertion is immediate for \(0\le t\le2\).
Let \(t>2\).
On \(0\le s\le t/2\), we have \(t-s\ge t/2\), and hence
\begin{equation*}
\begin{aligned}
&\int_0^{t/2}
\langle 2s-t\rangle^{-(1-\theta+\beta)}
\langle t-s\rangle^{-\theta}
\,\rmd s \\
&\qquad\le
C\langle t\rangle^{-\theta}
\int_0^{t/2}
\langle 2s-t\rangle^{-(1-\theta+\beta)}
\,\rmd s
\le C\langle t\rangle^{-\beta}.
\end{aligned}
\end{equation*}
On \(t/2\le s\le t\), putting \(u=t-s\), we estimate
\begin{equation*}
\int_0^{t/2}
\langle t-2u\rangle^{-(1-\theta+\beta)}
\langle u\rangle^{-\theta}
\,\rmd u.
\end{equation*}
For \(0\le u\le t/4\), \(t-2u\ge t/2\), whence
\begin{equation*}
\int_0^{t/4}
\langle t-2u\rangle^{-(1-\theta+\beta)}
\langle u\rangle^{-\theta}
\,\rmd u
\le C\langle t\rangle^{-(1-\theta+\beta)}
\langle t\rangle^{1-\theta}
\le C\langle t\rangle^{-\beta}.
\end{equation*}
For \(t/4\le u\le t/2\), \(\langle u\rangle\ge C^{-1}\langle t\rangle\), and therefore
\begin{equation*}
\begin{aligned}
&\int_{t/4}^{t/2}
\langle t-2u\rangle^{-(1-\theta+\beta)}
\langle u\rangle^{-\theta}
\,\rmd u \\
&\qquad\le
C\langle t\rangle^{-\theta}
\int_{t/4}^{t/2}
\langle t-2u\rangle^{-(1-\theta+\beta)}
\,\rmd u
\le C\langle t\rangle^{-\beta}.
\end{aligned}
\end{equation*}
This proves \eqref{eq:weak-wave-convolution}.

Consequently,
\begin{align*}
&\langle t+|x|\rangle^\theta\langle t-|x|\rangle^\beta
\times
\text{(left-hand side of \eqref{eq:weak-wave-duhamel-1d})}\\[1mm]
&\le C
\langle t+|x|\rangle^\theta\langle t-|x|\rangle^\beta
\langle t\rangle^{-\theta}
\int_0^t
\langle t-s\rangle^{-\theta}
\langle 2s-t\rangle^{-(1-\theta+\beta)}
\,\rmd s.
\end{align*}
By \eqref{eq:weak-wave-convolution},
\begin{equation*}
\int_0^t
\langle t-s\rangle^{-\theta}
\langle 2s-t\rangle^{-(1-\theta+\beta)}
\,\rmd s
\le C\langle t\rangle^{-\beta}.
\end{equation*}
Since \(|x|\le t/3\), we have
\begin{equation*}
\langle t+|x|\rangle\le C\langle t\rangle,
\qquad
\langle t-|x|\rangle\le C\langle t\rangle.
\end{equation*}
Therefore the last expression is bounded by a constant.
This proves \eqref{eq:weak-wave-duhamel-1d} in the case \(|x|\le (t+|x|)/4\).

It remains to consider \(|x|>(t+|x|)/4\).
Then \(|x|\) is comparable with \(t+|x|\).
We use the change of variables
\begin{equation*}
\mu=s+\lambda,
\qquad
\nu=s-\lambda.
\end{equation*}
The Jacobian is bounded by an absolute constant, and the cone condition
\begin{equation*}
||x|-(t-s)|\le \lambda\le |x|+t-s
\end{equation*}
implies
\begin{equation*}
|t-|x||\le \mu\le t+|x|,
\qquad
\nu\le t-|x|.
\end{equation*}
Since \(|x|>(t+|x|)/4\), we have \(|x|^{-1}\le C\langle t+|x|\rangle^{-1}\).
Moreover, \(\lambda=(\mu-\nu)/2\) and \(\lambda\le\langle\lambda\rangle\), and hence
\begin{equation*}
\frac{1}{|x|}\lambda\langle\lambda\rangle^{-(1+\theta)}
\le C\langle t+|x|\rangle^{-1}\langle\mu-\nu\rangle^{-\theta}.
\end{equation*}
Therefore the left-hand side of \eqref{eq:weak-wave-duhamel-1d} is bounded by
\begin{align*}
& C\langle t+|x|\rangle^{-1}
\int_{|t-|x||}^{t+|x|}
\int_{-\infty}^{t-|x|}
\langle\mu\rangle^{-\theta}
\langle\nu\rangle^{-(1-\theta+\beta)}
\langle\mu-\nu\rangle^{-\theta}\,\rmd\nu\,\rmd\mu.
\end{align*}
For each fixed \(\mu\in[|t-|x||,t+|x|]\), the inner integral is estimated by
\begin{equation*}
\int_{-\infty}^{t-|x|}
\langle\nu\rangle^{-(1-\theta+\beta)}
\langle\mu-\nu\rangle^{-\theta}\,\rmd\nu
\le C\langle t-|x|\rangle^{-\beta}.
\end{equation*}
This follows again from the one-dimensional convolution estimate, using \(0<\beta<\theta<1\).
Consequently,
\begin{align*}
\text{left-hand side of }\eqref{eq:weak-wave-duhamel-1d}
&\le C\langle t+|x|\rangle^{-1}
\langle t-|x|\rangle^{-\beta}
\int_{|t-|x||}^{t+|x|}\langle\mu\rangle^{-\theta}\,\rmd\mu \\
&\le C\langle t+|x|\rangle^{-\theta}
\langle t-|x|\rangle^{-\beta},
\end{align*}
where the last step uses \(0<\theta<1\).
This proves \eqref{eq:weak-wave-duhamel-1d}.
The case \(x=0\) is obtained by the same argument from the direct formula
\begin{equation*}
\LW[F](t,0)=\int_0^t (t-s)M_{t-s}F(s,\cdot)(0)\,\rmd s,
\end{equation*}
which is exactly the first case above.
Combining \eqref{eq:weak-wave-duhamel-reduced} and \eqref{eq:weak-wave-duhamel-1d} proves \eqref{eq:wave-inhom-weak}.
\end{proof}

For the damped wave equation, we use the following homogeneous estimate.

\begin{proposition}[\cite{Narazaki-Nishihara,Georgiev-Kita}]\label{prop:damped-hom}
Let \(0<\alpha<3/2\).
Then
\begin{equation}\label{eq:damped-hom-weak}
\|\KD[u_0,u_1]\|_{X_1^\alpha}
\le C\|(u_0,u_1)\|_{Y_1^\alpha}.
\end{equation}
If \(\alpha>3/2\), then
\begin{equation}\label{eq:damped-hom-strong}
\|\KD[u_0,u_1]\|_{X_1^{3/2}}
\le C\|(u_0,u_1)\|_{Y_1^\alpha}.
\end{equation}
\end{proposition}

\begin{remark}
The preceding proposition is a standard consequence of the explicit representation \eqref{eq:D-def}--\eqref{eq:S-def} and the pointwise estimates for \(D(t)\) and \(\partial_tD(t)\); see, for instance, \cite{Narazaki-Nishihara,Georgiev-Kita}.
At the endpoint, under only \(Y_1^{3/2}\)-decay, one generally has the logarithmic estimate
\begin{equation*}
|\KD[u_0,u_1](t,x)|
\le
C\langle t+|x|^2\rangle^{-3/2}
\log(e+t+|x|^2)
\|(u_0,u_1)\|_{Y_1^{3/2}}.
\end{equation*}
Thus \(Y_1^{3/2}\)-data are naturally covered by the weaker estimates \eqref{eq:damped-hom-weak} for every \(0<\alpha<3/2\), while the logarithmic endpoint estimate itself is not used in the proof of the main theorem.
In the strong case, where the conclusion is \(X_1^{3/2}\)-decay without logarithmic loss, we assume slightly faster spatial decay of the initial data.
\end{remark}

\section{Linear estimates}\label{sec:proof-damped}

This section proves Theorem \ref{thm:damped-duhamel}.
We write
\begin{equation}\label{eq:LD-split}
\LD[F]=L_{D,W}[F]+L_{D,S}[F],
\end{equation}
where
\begin{equation*}
L_{D,W}[F](t,x):=\int_0^t e^{-(t-\tau)/2}W(t-\tau)F(\tau,\cdot)(x)\,\rmd\tau,
\end{equation*}
\begin{equation*}
L_{D,S}[F](t,x):=\int_0^t S(t-\tau)F(\tau,\cdot)(x)\,\rmd\tau.
\end{equation*}
The first term is easier to estimate because of the exponential factor.
The main point is the estimate of the second term, where the kernel has a heat-like part but the source term contains the wave-cone weight.

We use the elementary identities
\begin{equation*}
I_1'(r)=I_0(r)-\frac1rI_1(r),\qquad I_0(r)-\frac2rI_1(r)=I_2(r),
\end{equation*}
and the standard estimates
\begin{equation*}
I_\nu(r)\le C
\left\{
\begin{aligned}
&r^\nu, &&0\le r\le \frac12,\\
&r^{-1/2}e^r, &&r>\frac12.
\end{aligned}
\right.
\end{equation*}

\begin{lemma}\label{lem:kernel-bound}
Set, for \(0\le \rho\le t\),
\begin{equation*}
\mathcal K(t,\rho):=e^{-t/2}\frac{I_1\left(\frac{\sqrt{t^2-\rho^2}}{2}\right)}{\sqrt{t^2-\rho^2}}.
\end{equation*}
Then
\begin{equation}\label{eq:kernel-crude}
\mathcal K(t,\rho)
\le C\frac{e^{-\rho^2/(4t)}}{\langle t^2-\rho^2\rangle^{3/4}}
\le C.
\end{equation}
Moreover, if \(t>\sqrt2\), then
\begin{equation}\label{eq:kernel-split}
\mathcal K(t,\rho)
\le C
\left\{
\begin{aligned}
&e^{-t/2}, &&\sqrt{t^2-\rho^2}\le 1,\\
&e^{-t/8}, &&\sqrt{t^2-\rho^2}>1,\quad t/\sqrt2\le \rho,\\
&t^{-3/2}e^{-\rho^2/(4t)}, &&\sqrt{t^2-\rho^2}>1,\quad t/\sqrt2>\rho.
\end{aligned}
\right.
\end{equation}
\end{lemma}

\begin{proof}
Let \(r=\sqrt{t^2-\rho^2}/2\).
If \(r\le 1/2\), then \(I_1(r)\le Cr\), and hence \(\mathcal K(t,\rho)\le Ce^{-t/2}\).
If \(r>1/2\), then \(I_1(r)\le Cr^{-1/2}e^r\), and therefore
\begin{equation*}
\mathcal K(t,\rho)
\le C(t^2-\rho^2)^{-3/4}
\exp\left(-\frac{t-\sqrt{t^2-\rho^2}}2\right).
\end{equation*}
Since
\begin{equation*}
t-\sqrt{t^2-\rho^2}=\frac{\rho^2}{t+\sqrt{t^2-\rho^2}}\ge \frac{\rho^2}{2t},
\end{equation*}
we get \eqref{eq:kernel-crude}.
If \(\rho\ge t/\sqrt2\), then \(t-\sqrt{t^2-\rho^2}\ge c t\), and hence the exponential factor is bounded by \(e^{-ct}\), which gives the second line of \eqref{eq:kernel-split}.
If \(\rho<t/\sqrt2\), then \(t^2-\rho^2\ge t^2/2\), and the last line of \eqref{eq:kernel-split} follows.
\end{proof}

\begin{lemma}\label{lem:LDW-part}
Let \(\alpha\in(0,3/2]\), \(A\ge \alpha\), and \(\gamma\ge 0\).
Then
\begin{equation*}
\|L_{D,W}[F]\|_{X_1^\alpha}
\le C
\left\|\langle t+|x|^2\rangle^A\langle t-|x|\rangle^\gamma F(t,x)\right\|_{L^\infty_{t,x}}.
\end{equation*}
\end{lemma}

\begin{proof}
For every \(\kappa\ge0\), by Lemma 2.1 of \cite{Georgiev-Kita}, the spherical mean estimate and the exponential factor give
\begin{equation} \label{eq:LDW-spatial-kappa}
\left|e^{-s/2}W(s)f(x)\right|
\le C e^{-s/8}\langle x\rangle^{-\kappa}\|f\|_{\infty,\kappa}.
\end{equation}
Applying \eqref{eq:LDW-spatial-kappa} with \(s=t-\tau\), \(f=F(\tau,\cdot)\), and \(\kappa=2\alpha\), we obtain
\begin{equation*}
\left|e^{-(t-\tau)/2}W(t-\tau)F(\tau,\cdot)(x)\right|
\le C e^{-c(t-\tau)}\langle x\rangle^{-2\alpha}
\|F(\tau,\cdot)\|_{\infty,2\alpha}.
\end{equation*}
Since \(A\ge\alpha\), \(\gamma\ge0\) and \(\langle \tau+|y|^2\rangle^A\langle \tau-|y|\rangle^\gamma \ge C\langle y\rangle^{2\alpha}\), we have
\begin{equation*}
\|F(\tau,\cdot)\|_{\infty,2\alpha}
\le C M,
\end{equation*}
where
\begin{equation*}
    M:=
\left\|
\langle s+|y|^2\rangle^A
\langle s-|y|\rangle^\gamma
F(s,y)
\right\|_{L^\infty_{s,y}}.
\end{equation*}
Therefore
\begin{equation*}
|L_{D,W}[F](t,x)|
\le C\langle x\rangle^{-2\alpha}M
\int_0^t e^{-c(t-\tau)}\,\rmd\tau
\le C\langle x\rangle^{-2\alpha}M.
\end{equation*}
This yields the spatial decay estimate
\begin{equation*}
|L_{D,W}[F](t,x)|\le CM\langle x\rangle^{-2\alpha}.
\end{equation*}

On the other hand, using only the time weight in the source and the exponential integrability in \(t-\tau\), we have
\begin{align*}
|L_{D,W}[F](t,x)|
&\le C\int_0^t e^{-c(t-\tau)}\|F(\tau,\cdot)\|_{L^\infty}\,\rmd\tau \\
&\le C\int_0^t e^{-c(t-\tau)}\langle \tau\rangle^{-A}\,\rmd\tau
\left\|\langle s+|y|^2\rangle^A\langle s-|y|\rangle^\gamma F(s,y)\right\|_{L^\infty_{s,y}} \\
&\le C\langle t\rangle^{-A}
\left\|\langle s+|y|^2\rangle^A\langle s-|y|\rangle^\gamma F(s,y)\right\|_{L^\infty_{s,y}} \\
&\le C\langle t\rangle^{-\alpha}
\left\|\langle s+|y|^2\rangle^A\langle s-|y|\rangle^\gamma F(s,y)\right\|_{L^\infty_{s,y}}.
\end{align*}
Here the factor \(t-\tau\) in \(W(t-\tau)\) is absorbed by the exponential decay.
Hence the left-hand side is bounded by the minimum of the two right-hand sides.
Since
\begin{equation*}
\min\{\langle t\rangle^{-\alpha},\langle x\rangle^{-2\alpha}\}
\le C\langle t+|x|^2\rangle^{-\alpha},
\end{equation*}
the claim follows.
\end{proof}

We now estimate the \(S(t)\)-part.
Put
\begin{equation*}
M_{A,\gamma}(F):=
\left\|\langle t+|x|^2\rangle^A\langle t-|x|\rangle^\gamma F(t,x)\right\|_{L^\infty_{t,x}}.
\end{equation*}
From the explicit formula for \(S(t)\), Lemma \ref{lem:kernel-bound} gives
\begin{equation}\label{eq:LS-basic-bound}
|L_{D,S}[F](t,x)|
\le C M_{A,\gamma}(F)\mathcal I_{A,\gamma}(t,x),
\end{equation}
where
\begin{equation*}
\mathcal I_{A,\gamma}(t,x):=\int_0^t\int_0^{t-\tau}\int_{\mathbb S^2}
\mathcal K(t-\tau,\rho)\rho^2
\langle \tau+|x+\rho\omega|^2\rangle^{-A}
\langle \tau-|x+\rho\omega|\rangle^{-\gamma}
\,\rmd\sigma(\omega)\,\rmd\rho\,\rmd\tau.
\end{equation*}
The following proposition is the estimate needed for both the strong and weak versions of the damped wave Duhamel estimate.

\begin{proposition}\label{prop:S-part}
Let \(0<\alpha\le 3/2\), \(A>\alpha\), and \(\gamma>1\).
Then
\begin{equation}\label{eq:I-unified}
\mathcal I_{A,\gamma}(t,x)
\le C\langle t+|x|^2\rangle^{-\alpha}
\end{equation}
for all \(t\ge0\) and \(x\in\R^3\).
Consequently,
\begin{equation*}
|L_{D,S}[F](t,x)|
\le C\langle t+|x|^2\rangle^{-\alpha}M_{A,\gamma}(F).
\end{equation*}
\end{proposition}

\begin{proof}
The proof is divided into the parabolic region \(|x|^2\le t\) and the exterior region \(t<|x|^2\).
This avoids an unnecessary loss in the weak range \(0<\alpha<3/2\).
Constants may depend on \(\alpha,A,\gamma\), but not on \(t\), \(x\), or \(F\).

If \(0\le t\le 2\sqrt2\), then \(0\le \rho\le t\le 2\sqrt2\).
By \eqref{eq:kernel-crude},
\begin{equation*}
\mathcal I_{A,\gamma}(t,x)
\le C\int_0^t\int_0^{t-\tau}\int_{\mathbb S^2}\rho^2
\langle \tau+|x+\rho\omega|^2\rangle^{-A}\,\rmd\sigma(\omega)\,\rmd\rho\,\rmd\tau.
\end{equation*}
For bounded \(|x|\), the desired estimate is immediate, since \(\langle t+|x|^2\rangle\) is bounded.
If \(|x|\ge 4\sqrt2\), then \(|x+\rho\omega|\ge |x|/2\), and hence the last integral is bounded by \(C\langle x\rangle^{-2A}\).
Since \(A>\alpha\), this is bounded by \(C\langle t+|x|^2\rangle^{-\alpha}\) for \(0\le t\le2\sqrt2\).
Henceforth we assume \(t>2\sqrt2\).

We decompose the integration region according to the kernel estimate \eqref{eq:kernel-split}.
Let
\begin{equation*}
\Omega_0:=\{(\tau,\rho):t-\sqrt2\le \tau\le t,\ 0\le \rho\le t-\tau\},
\end{equation*}
\begin{equation*}
\Omega_1:=\{(\tau,\rho):0\le\tau\le t-\sqrt2,
\sqrt{(t-\tau)^2-\rho^2}\le 1\},
\end{equation*}
\begin{equation*}
\Omega_2:=\{(\tau,\rho):0\le\tau\le t-\sqrt2,
\sqrt{(t-\tau)^2-\rho^2}>1,
(t-\tau)/\sqrt2\le \rho\le t-\tau\},
\end{equation*}
\begin{equation*}
\Omega_3:=\{(\tau,\rho):0\le\tau\le t-\sqrt2,
\sqrt{(t-\tau)^2-\rho^2}>1,
0\le \rho<(t-\tau)/\sqrt2\}.
\end{equation*}
Let \(\mathcal I_j\) be the contribution of \(\Omega_j\), \(j=0,1,2,3\).

We first estimate \(\mathcal I_0+\mathcal I_1+\mathcal I_2\).
In the parabolic region \(|x|^2\le t\), the term \(\mathcal I_0\) satisfies \(\tau\ge t-\sqrt2\), and the integration in \(\rho\) has bounded length.
Hence
\begin{equation*}
\mathcal I_0\le C\langle t\rangle^{-A}.
\end{equation*}
For \(\mathcal I_1+\mathcal I_2\), put \(\sigma=t-\tau\).
By \eqref{eq:kernel-split}, the kernel is bounded by \(C e^{-c\sigma}\), and the \(\rho\)-volume is bounded by \(C\sigma^3\).
Thus
\begin{equation*}
\mathcal I_1+\mathcal I_2
\le C\int_0^{t-\sqrt2} e^{-c(t-\tau)}(t-\tau)^3\langle \tau\rangle^{-A}\,\rmd\tau.
\end{equation*}
On \([0,t/2]\) the exponential factor gives \(C_N\langle t\rangle^{-N}\) for any fixed \(N>0\), while on \([t/2,t-\sqrt2]\) one has \(\langle \tau\rangle^{-A}\le C\langle t\rangle^{-A}\).
Hence
\begin{equation}\label{eq:I012-parabolic}
\mathcal I_0+\mathcal I_1+\mathcal I_2
\le C\langle t\rangle^{-A}+C_N\langle t\rangle^{-N}
\le C\langle t\rangle^{-\alpha}.
\end{equation}

In the exterior region \(t<|x|^2\), we use the spatial weight.
If \(|x|\) is bounded, the desired estimate is trivial.
We may therefore assume \(|x|>2\sqrt2\).
On \(\Omega_0\), \(\rho\le\sqrt2\), so \(|x+\rho\omega|\ge |x|/2\), and
\begin{equation*}
\mathcal I_0\le C\langle x\rangle^{-2A}.
\end{equation*}
For \(\mathcal I_1+\mathcal I_2\), split the radial interval into \(\rho>|x|/2\) and \(\rho\le |x|/2\).
If \(\rho>|x|/2\), then \(t-\tau\ge \rho>|x|/2\), and the exponential factor gives a contribution bounded by \(C_N\langle x\rangle^{-N}\).
If \(\rho\le |x|/2\), then \(|x+\rho\omega|\ge |x|/2\), and the source weight gives \(C\langle x\rangle^{-2A}\).
Therefore
\begin{equation}\label{eq:I012-exterior}
\mathcal I_0+\mathcal I_1+\mathcal I_2
\le C\langle x\rangle^{-2A}
\le C\langle x\rangle^{-2\alpha}.
\end{equation}

It remains to estimate \(\mathcal I_3\).
On \(\Omega_3\), since \(t-\tau\ge \sqrt{2}\), we have
\begin{equation*}
\frac{e^{-\rho^2/(4(t-\tau))}}{(t-\tau)^{3/2}}
\le C\langle\rho\rangle^{-3}.
\end{equation*}
We use the change of variables \(z=x+\rho\omega\).
Since \(\rho=|z-x|\), the Gaussian part of \eqref{eq:kernel-split} gives
\begin{equation}\label{eq:I3-z-form}
\mathcal I_3
\le C\int_0^{t-\sqrt2}\int_{|x-z|\le (t-\tau)/\sqrt2}
\frac{e^{-|x-z|^2/(4(t-\tau))}}{(t-\tau)^{3/2}}
\langle \tau+|z|^2\rangle^{-A}
\langle \tau-|z|\rangle^{-\gamma}\,\rmd z\,\rmd\tau.
\end{equation}

We use two elementary estimates.
The first one is used in the parabolic region.
If \(\mu\ge1\), then
\begin{equation}\label{eq:radial-cone-gauss}
\int_0^\mu\int_{\R^3}e^{-c|z|^2/\mu}
\langle \tau+|z|^2\rangle^{-A}
\langle \tau-|z|\rangle^{-\gamma}\,\rmd z\,\rmd\tau
\le C
\left\{
\begin{aligned}
&\mu^{3/2-A}, &&0<A<3/2,\\
&\log(e+\mu), &&A=3/2,\\
&1, &&A>3/2.
\end{aligned}
\right.
\end{equation}
By polar coordinates it suffices to estimate
\begin{equation*}
\int_0^\mu\int_0^\infty e^{-cr^2/\mu}r^2
\langle \tau+r^2\rangle^{-A}\langle \tau-r\rangle^{-\gamma}\,\rmd r\,\rmd\tau.
\end{equation*}
For fixed \(r\ge0\), we estimate the \(\tau\)-integral directly.
Since \(\tau\ge0\),
\begin{equation*}
\langle \tau+r^2\rangle^{-A}\le C\langle r^2\rangle^{-A}.
\end{equation*}
Moreover, since \(\gamma>1\),
\begin{equation*}
\int_0^\mu \langle \tau-r\rangle^{-\gamma}\,\rmd\tau
\le \int_0^\infty \langle \tau-r\rangle^{-\gamma}\,\rmd\tau
\le C.
\end{equation*}
Therefore
\begin{equation*}
\int_0^\mu
\langle \tau+r^2\rangle^{-A}\langle \tau-r\rangle^{-\gamma}\,\rmd\tau
\le C\langle r^2\rangle^{-A}.
\end{equation*}
Thus the left-hand side of \eqref{eq:radial-cone-gauss} is bounded by
\begin{equation*}
C\int_0^\infty e^{-cr^2/\mu}r^2\langle r^2\rangle^{-A}\,\rmd r.
\end{equation*}
Splitting the last integral at \(r=1\) and using the change of variables \(r=\sqrt\mu\,s\) on \(r\ge1\), we obtain the three alternatives in \eqref{eq:radial-cone-gauss}: the bound is \(C\mu^{3/2-A}\) if \(0<A<3/2\), \(C\log(e+\mu)\) if \(A=3/2\), and \(C\) if \(A>3/2\).

The second estimate is used in the exterior region.
For \(|x|\ge2\),
\begin{equation}\label{eq:tail-conv}
\int_{|x-z|\ge |x|/2}\langle x-z\rangle^{-3}\langle z\rangle^{-2A}\,\rmd z
\le C\langle x\rangle^{-2\alpha}.
\end{equation}
To prove this, split the integral into the three regions \(|z|\le |x|/2\), \(|x|/2<|z|<2|x|\), and \(|z|\ge2|x|\).
In the first region,
\begin{equation*}
\int_{|z|\le |x|/2}\langle x-z\rangle^{-3}\langle z\rangle^{-2A}\,\rmd z
\le C\langle x\rangle^{-3}\int_{|z|\le |x|/2}\langle z\rangle^{-2A}\,\rmd z.
\end{equation*}
This is bounded by \(C\langle x\rangle^{-2A}\) if \(A<3/2\), by \(C\langle x\rangle^{-3}\log(e+|x|)\) if \(A=3/2\), and by \(C\langle x\rangle^{-3}\) if \(A>3/2\).
Each of these bounds is at most \(C\langle x\rangle^{-2\alpha}\), since \(A>\alpha\) and \(\alpha\le3/2\).
In the middle region, \(\langle z\rangle^{-2A}\le C\langle x\rangle^{-2A}\), and
\begin{equation*}
\int_{\substack{|x|/2<|z|<2|x|\\ |x-z|\ge |x|/2}}\langle x-z\rangle^{-3}\,\rmd z\le C.
\end{equation*}
The contribution is therefore bounded by \(C\langle x\rangle^{-2A}\).
In the last region, \(\langle x-z\rangle\ge c\langle z\rangle\), and the contribution is bounded by
\begin{equation*}
C\int_{|z|\ge2|x|}\langle z\rangle^{-3-2A}\,\rmd z
\le C\langle x\rangle^{-2A}.
\end{equation*}
This proves \eqref{eq:tail-conv}.

We now prove the parabolic estimate for \(\mathcal I_3\).
Assume \(|x|^2\le t\).
Split \(\mathcal I_3=\mathcal I_3^{\rm e}+\mathcal I_3^{\rm l}\), where \(\mathcal I_3^{\rm e}\) corresponds to \(0\le\tau\le t/2\), and \(\mathcal I_3^{\rm l}\) corresponds to \(t/2\le\tau\le t-\sqrt2\).

For \(\mathcal I_3^{\rm e}\), we have \(t-\tau\ge t/2\).
Since \(|x|^2\le t\),
\begin{equation*}
\exp\left(-\frac{|x-z|^2}{4(t-\tau)}\right)
\le C\exp\left(-c\frac{|z|^2}{t}\right).
\end{equation*}
Thus \eqref{eq:I3-z-form} and \eqref{eq:radial-cone-gauss} with \(\mu=t\) give
\begin{equation*}
\mathcal I_3^{\rm e}
\le C t^{-3/2}
\int_0^{t/2}\int_{\R^3}e^{-c|z|^2/t}
\langle \tau+|z|^2\rangle^{-A}
\langle \tau-|z|\rangle^{-\gamma}\,\rmd z\,\rmd\tau
\le C\langle t\rangle^{-\alpha}.
\end{equation*}
Indeed, if \(A<3/2\), the bound is \(C t^{-A}\), and \(A>\alpha\).
If \(A=3/2\), then necessarily \(\alpha<3/2\), and
\begin{equation*}
t^{-3/2}\log(e+t)\le C\langle t\rangle^{-\alpha}.
\end{equation*}
If \(A>3/2\), the bound is \(C t^{-3/2}\), which is enough since \(\alpha\le3/2\).

For \(\mathcal I_3^{\rm l}\), we have \(\tau\ge t/2\), hence
\begin{equation*}
\langle \tau+|z|^2\rangle^{-A}\le C\langle t\rangle^{-A}.
\end{equation*}
It remains to show that the heat kernel integrated against the cone weight is bounded by a constant independent of \(t\).
We split the \(z\)-integral into \(|z|\le \tau/2\) and \(|z|>\tau/2\).
On \(|z|\le \tau/2\), we have \(\langle \tau-|z|\rangle^{-\gamma}\le C\langle \tau\rangle^{-\gamma}\), and the heat kernel has bounded mass.
Hence this part contributes at most
\begin{equation*}
C\langle t\rangle^{-A}\int_{t/2}^{t}\langle \tau\rangle^{-\gamma}\,\rmd\tau
\le C\langle t\rangle^{-A}.
\end{equation*}
On \(|z|>\tau/2\), use \(|x|^2\le t\) and \(\tau\ge t/2\).
For sufficiently large \(t\), this gives \(|x-z|\ge c\tau\); the case of bounded \(t\) is immediate.
Hence, for any \(N>0\),
\begin{equation*}
\int_{|z|>\tau/2}\frac{e^{-|x-z|^2/(4(t-\tau))}}{(t-\tau)^{3/2}}\,\rmd z
\le C_N\left(1+\frac{\tau^2}{t-\tau}\right)^{-N}.
\end{equation*}
Putting \(\sigma=t-\tau\) and changing variables \(y=(z-x)/\sqrt\sigma\), the region under consideration is contained in \(|y|\ge c\tau/\sqrt\sigma\) for large \(t\), and the Gaussian tail satisfies
\begin{equation*}
\int_{|y|\ge c\tau/\sqrt\sigma}e^{-|y|^2/4}\,\rmd y
\le C_N\left(1+\frac{\tau^2}{\sigma}\right)^{-N}.
\end{equation*}
Integrating over \(\tau\in[t/2,t-\sqrt2]\) and taking \(N\) large, this gives another contribution bounded by \(C\langle t\rangle^{-A}\).
Since \(A>\alpha\), we obtain
\begin{equation*}
\mathcal I_3\le C\langle t\rangle^{-\alpha}
\qquad (|x|^2\le t).
\end{equation*}
Combining this with \eqref{eq:I012-parabolic} proves \eqref{eq:I-unified} in the parabolic region.

It remains to prove the exterior estimate.
Assume \(t<|x|^2\) and \(|x|>1\).
By \eqref{eq:I012-exterior}, only \(\mathcal I_3\) remains.
In \eqref{eq:I3-z-form}, split the spatial region into
\begin{equation*}
U^e:=\{z:\ |x-z|\ge |x|/2\},
\qquad
U^i:=\{z:\ |x-z|< |x|/2\}.
\end{equation*}
Let \(\mathcal I_3^e\) and \(\mathcal I_3^i\) denote the corresponding contributions.

For \(\mathcal I_3^e\), use
\begin{equation*}
\sup_{\sigma>0}\sigma^{-3/2}\exp\left(-\frac{|x-z|^2}{4\sigma}\right)
\le C\langle x-z\rangle^{-3}.
\end{equation*}
Moreover,
\begin{equation*}
\langle \tau+|z|^2\rangle^{-A}\le C\langle z\rangle^{-2A},
\qquad
\int_0^\infty\langle \tau-|z|\rangle^{-\gamma}\,\rmd\tau\le C.
\end{equation*}
Therefore \eqref{eq:tail-conv} yields
\begin{equation}\label{eq:I3e-exterior}
\mathcal I_3^e
\le C\int_{|x-z|\ge |x|/2}\langle x-z\rangle^{-3}\langle z\rangle^{-2A}\,\rmd z
\le C\langle x\rangle^{-2\alpha}.
\end{equation}

For \(\mathcal I_3^i\), we have \(|z|\ge |x|/2\), and hence
\begin{equation*}
\langle \tau+|z|^2\rangle^{-A}\le C\langle x\rangle^{-2A}.
\end{equation*}
We return to the variables \(\rho=|x-z|\) and \(\omega\).
On \(\Omega_3\), \(t-\tau\ge\sqrt2\), and therefore
\begin{equation*}
\frac{e^{-\rho^2/(4(t-\tau))}}{(t-\tau)^{3/2}}
\le C\langle\rho\rangle^{-3}.
\end{equation*}
Using the spherical identity
\begin{equation*}
\frac{\rho}{4\pi}\int_{\mathbb S^2}h(|x+\rho\omega|)\,\rmd\sigma(\omega)
=\frac{1}{2|x|}\int_{\bigl||x|-\rho\bigr|}^{|x|+\rho}\lambda h(\lambda)\,\rmd\lambda,
\end{equation*}
we get
\begin{align*}
\mathcal I_3^i
&\le C\langle x\rangle^{-2A}
\int_0^{|x|/2}\langle\rho\rangle^{-3}\rho^2
\int_{\mathbb S^2}\int_0^\infty
\langle \tau-|x+\rho\omega|\rangle^{-\gamma}\,\rmd\tau\,\rmd\sigma(\omega)\,\rmd\rho \\
&\le C\langle x\rangle^{-2A}
\int_0^{|x|/2}\langle\rho\rangle^{-2}
\left\{\int_{\bigl||x|-\rho\bigr|}^{|x|+\rho}\frac{\lambda}{|x|}
\left(\int_0^\infty\langle\tau-\lambda\rangle^{-\gamma}\,\rmd\tau\right)\,\rmd\lambda\right\}\,\rmd\rho.
\end{align*}
Since \(\gamma>1\), the inner \(\tau\)-integral is uniformly bounded.
Also, for \(0\le\rho\le |x|/2\),
\begin{equation*}
\int_{\bigl||x|-\rho\bigr|}^{|x|+\rho}\frac{\lambda}{|x|}\,\rmd\lambda\le C\rho.
\end{equation*}
Thus
\begin{equation}\label{eq:I3i-exterior}
\mathcal I_3^i
\le C\langle x\rangle^{-2A}\int_0^{|x|/2}\langle\rho\rangle^{-2}\rho\,\rmd\rho
\le C\langle x\rangle^{-2A}\log(e+|x|)
\le C\langle x\rangle^{-2\alpha}.
\end{equation}
Combining \eqref{eq:I3e-exterior} and \eqref{eq:I3i-exterior}, we obtain
\begin{equation*}
\mathcal I_3\le C\langle x\rangle^{-2\alpha}
\qquad (t<|x|^2).
\end{equation*}
Together with \eqref{eq:I012-exterior}, this proves \eqref{eq:I-unified} in the exterior region.
This completes the proof of the proposition.
\end{proof}

\begin{proof}[Proof of Theorem \ref{thm:damped-duhamel}]
Let
\begin{equation*}
M_{A,\gamma}(F)
:=
\left\|\langle t+|x|^2\rangle^A\langle t-|x|\rangle^\gamma F(t,x)\right\|_{L^\infty_{t,x}}.
\end{equation*}
By the decomposition \eqref{eq:LD-split},
\begin{equation*}
\LD[F]=L_{D,W}[F]+L_{D,S}[F].
\end{equation*}
Hence, by the triangle inequality in \(X_1^\alpha\),
\begin{equation*}
\|\LD[F]\|_{X_1^\alpha}
\le
\|L_{D,W}[F]\|_{X_1^\alpha}
+
\|L_{D,S}[F]\|_{X_1^\alpha}.
\end{equation*}

The first term is estimated by Lemma \ref{lem:LDW-part}:
\begin{equation*}
\|L_{D,W}[F]\|_{X_1^\alpha}\le C M_{A,\gamma}(F).
\end{equation*}
For the second term, \eqref{eq:LS-basic-bound} together with Proposition \ref{prop:S-part} yields
\begin{equation*}
\|L_{D,S}[F]\|_{X_1^\alpha}\le C M_{A,\gamma}(F).
\end{equation*}
Therefore
\begin{equation*}
\|\LD[F]\|_{X_1^\alpha}\le C M_{A,\gamma}(F),
\end{equation*}
which is exactly \eqref{eq:damped-duhamel}.
This completes the proof.
\end{proof}

\section{Global existence for Nakao's problem}\label{sec:global}

We prove Theorem \ref{thm:global-parameter} and then derive Corollaries \ref{cor:strong-strong}--\ref{cor:strong-weak}.
The proof is a contraction argument, but the choice of auxiliary exponents is important.

\begin{proof}[Proof of Theorem \ref{thm:global-parameter}]
By \eqref{eq:param-cond}, we can choose \(A\) and \(\gamma\) such that
\begin{equation}\label{eq:A-gamma-choice}
\alpha<A<\frac{\theta p}{2},
\qquad
1<\gamma<\beta p.
\end{equation}
If \(\theta=1\), we also choose \(m>2\) so close to \(2\) that
\begin{equation}\label{eq:m-choice}
1+\beta+\frac m2<\alpha q.
\end{equation}
This is possible since \(\alpha q>2+\beta\) when \(\theta=1\).

For \((u,v)\in X_1^\alpha\times X_2^{\theta,\beta}\), define
\begin{equation*}
\Phi(u,v):=\left(\KD[u_0,u_1]+\LD[|v|^p],\ \KW[v_0,v_1]+\LW[|u|^q]\right).
\end{equation*}
We first estimate the nonlinear terms.
If \(v\in X_2^{\theta,\beta}\), then
\begin{equation*}
|v(t,x)|^p\le \|v\|_{X_2^{\theta,\beta}}^p
\langle t+|x|\rangle^{-\theta p}\langle t-|x|\rangle^{-\beta p}.
\end{equation*}
Using \(\langle t+|x|^2\rangle\le C\langle t+|x|\rangle^2\), \eqref{eq:A-gamma-choice} gives
\begin{equation}\label{eq:v-to-u-source}
\left\|\langle t+|x|^2\rangle^A\langle t-|x|\rangle^\gamma |v|^p\right\|_{L^\infty_{t,x}}
\le C\|v\|_{X_2^{\theta,\beta}}^p.
\end{equation}
Hence Theorem \ref{thm:damped-duhamel} yields
\begin{equation}\label{eq:u-nonlinear-bound}
\|\LD[|v|^p]\|_{X_1^\alpha}
\le C\|v\|_{X_2^{\theta,\beta}}^p.
\end{equation}

Next, if \(u\in X_1^\alpha\), then
\begin{equation*}
|u(t,x)|^q\le \|u\|_{X_1^\alpha}^q\langle t+|x|^2\rangle^{-\alpha q}.
\end{equation*}
When \(0<\theta<1\), the weak wave estimate \eqref{eq:wave-inhom-weak} and \eqref{eq:param-cond} give
\begin{equation}\label{eq:u-to-v-source-weak}
\left\|\langle t+|x|\rangle^\theta\langle t-|x|\rangle^{1-\theta+\beta}\langle x\rangle^{1+\theta}|u|^q\right\|_{L^\infty_{t,x}}
\le C\|u\|_{X_1^\alpha}^q.
\end{equation}
Indeed,
\begin{equation*}
\langle t+|x|\rangle^\theta\langle t-|x|\rangle^{1-\theta+\beta}\langle x\rangle^{1+\theta}
\le C\langle t+|x|^2\rangle^{1+\beta+(1+\theta)/2}.
\end{equation*}
Therefore
\begin{equation}\label{eq:v-nonlinear-bound-weak}
\|\LW[|u|^q]\|_{X_2^{\theta,\beta}}
\le C\|u\|_{X_1^\alpha}^q.
\end{equation}
When \(\theta=1\), we instead use \eqref{eq:wave-inhom-strong} with the \(m\) chosen in \eqref{eq:m-choice}.
Since \(\beta<1\), this \(m\) may also be chosen so that \(\beta<m-1\), as required in Proposition \ref{prop:wave-estimates}.
The estimate
\begin{equation*}
\left\|\langle t+|x|\rangle\langle t-|x|\rangle^\beta\langle x\rangle^m|u|^q\right\|_{L^\infty_{t,x}}
\le C\|u\|_{X_1^\alpha}^q
\end{equation*}
follows from \(1+\beta+m/2<\alpha q\).
Thus \eqref{eq:v-nonlinear-bound-weak} remains valid with \(\theta=1\), after replacing the weak wave estimate by the strong one.

The homogeneous estimates \eqref{eq:damped-hom-weak}, \eqref{eq:damped-hom-strong}, \eqref{eq:wave-hom-strong}, and \eqref{eq:wave-hom-weak} imply
\begin{equation}\label{eq:hom-bound-global}
\|\KD[u_0,u_1]\|_{X_1^\alpha}+\|\KW[v_0,v_1]\|_{X_2^{\theta,\beta}}
\le C\varepsilon,
\end{equation}
where
\begin{equation*}
\varepsilon:=\mathcal N_0.
\end{equation*}
Combining \eqref{eq:u-nonlinear-bound}, \eqref{eq:v-nonlinear-bound-weak}, and \eqref{eq:hom-bound-global}, we obtain
\begin{equation}\label{eq:Phi-bound}
\|\Phi(u,v)\|_{X_1^\alpha\times X_2^{\theta,\beta}}
\le C\varepsilon+C\|u\|_{X_1^\alpha}^q+C\|v\|_{X_2^{\theta,\beta}}^p.
\end{equation}
We also need the corresponding Lipschitz estimates.
Let \((u,v)\) and \((\tilde u,\tilde v)\) belong to the ball of radius \(\varrho\) in \(X_1^\alpha\times X_2^{\theta,\beta}\).
By the mean value inequality,
\begin{equation*}
\bigl||v|^p-|\tilde v|^p\bigr|
\le C\left(|v|^{p-1}+|\tilde v|^{p-1}\right)|v-\tilde v|,
\end{equation*}
and the same weight computation as in \eqref{eq:v-to-u-source} gives
\begin{equation}\label{eq:v-lip-source}
\left\|\langle t+|x|^2\rangle^A\langle t-|x|\rangle^\gamma
\bigl(|v|^p-|\tilde v|^p\bigr)\right\|_{L^\infty_{t,x}}
\le C\varrho^{p-1}\|v-\tilde v\|_{X_2^{\theta,\beta}}.
\end{equation}
Similarly,
\begin{equation*}
\bigl||u|^q-|\tilde u|^q\bigr|
\le C\left(|u|^{q-1}+|\tilde u|^{q-1}\right)|u-\tilde u|.
\end{equation*}
If \(0<\theta<1\), then the same computation as in \eqref{eq:u-to-v-source-weak} yields
\begin{equation*}
\left\|\langle t+|x|\rangle^\theta\langle t-|x|\rangle^{1-\theta+\beta}
\langle x\rangle^{1+\theta}
\bigl(|u|^q-|\tilde u|^q\bigr)\right\|_{L^\infty_{t,x}}
\le C\varrho^{q-1}\|u-\tilde u\|_{X_1^\alpha}.
\end{equation*}
If \(\theta=1\), the analogous estimate with the strong wave weight
\begin{equation}\label{eq:u-lip-source-strong}
\left\|\langle t+|x|\rangle\langle t-|x|\rangle^\beta
\langle x\rangle^m
\bigl(|u|^q-|\tilde u|^q\bigr)\right\|_{L^\infty_{t,x}}
\le C\varrho^{q-1}\|u-\tilde u\|_{X_1^\alpha}
\end{equation}
follows from \(1+\beta+m/2<\alpha q\).
Applying Theorem \ref{thm:damped-duhamel} and the wave estimates to \eqref{eq:v-lip-source}--\eqref{eq:u-lip-source-strong}, we obtain
\begin{equation}\label{eq:Phi-Lip}
\|\Phi(u,v)-\Phi(\tilde u,\tilde v)\|_{X_1^\alpha\times X_2^{\theta,\beta}}
\le C\left(\varrho^{p-1}+\varrho^{q-1}\right)
\|(u-\tilde u,v-\tilde v)\|_{X_1^\alpha\times X_2^{\theta,\beta}}.
\end{equation}
Let \(C_0\) be a constant such that both \eqref{eq:Phi-bound} and \eqref{eq:Phi-Lip} hold.
Choose a fixed \(K\ge 2C_0\) and set \(\varrho=K\varepsilon\).
If \(\varepsilon>0\) is sufficiently small, depending only on the displayed exponents and on \(K\), then
\begin{equation*}
C_0\varepsilon+C_0\varrho^p+C_0\varrho^q\le \varrho,
\qquad
C_0\left(\varrho^{p-1}+\varrho^{q-1}\right)<1.
\end{equation*}
Thus \(\Phi\) maps the closed ball of radius \(\varrho\) in \(X_1^\alpha\times X_2^{\theta,\beta}\) into itself and is a contraction there.
The contraction mapping theorem gives a unique fixed point in this ball.
This fixed point is the desired global mild solution, and \eqref{eq:global-estimate} follows with \(C=K\).
The same contraction estimate also gives uniqueness among all global mild solutions satisfying \eqref{eq:global-estimate} with this constant \(C\).
\end{proof}

\begin{proof}[Proof of Corollary \ref{cor:strong-strong}]
Set \(\alpha=3/2\) and \(\theta=1\).
The first condition in \eqref{eq:param-cond} is \(p>3\).
The condition \(\beta p>1\) requires \(\beta>1/p\).
The last condition is \((3/2)q>2+\beta\).
Such a \(\beta\in(1/p,1)\) exists if and only if
\begin{equation*}
\frac32 q>2+\frac1p,
\end{equation*}
which is equivalent to \eqref{eq:ss-cond}.
\end{proof}

\begin{proof}[Proof of Corollary \ref{cor:weak-strong}]
We use Theorem \ref{thm:global-parameter} with \(\theta=1\).
We need
\begin{equation*}
2\alpha<p,\qquad \beta>\frac1p,\qquad \alpha q>2+\beta.
\end{equation*}
Since \(1<p\le3\), choosing \(\alpha<p/2\) automatically gives \(\alpha<3/2\).
Hence such \(\alpha\in(0,3/2)\) and \(\beta\in(1/p,1)\) exist if and only if
\begin{equation*}
\frac p2 q>2+\frac1p,
\end{equation*}
which is equivalent to \eqref{eq:ws-cond}.
\end{proof}

\begin{proof}[Proof of Corollary \ref{cor:strong-weak}]
Set \(\alpha=3/2\).
We need
\begin{equation*}
\theta p>3,
\qquad
\beta p>1,
\qquad
\frac32 q>1+\beta+\frac{1+\theta}{2}.
\end{equation*}
The first two inequalities can be satisfied with \(0<\beta<\theta<1\) since \(p>3\).
Since the left inequality in \eqref{eq:sw-cond} is strict, we may take \(\theta>3/p\) and \(\beta>1/p\) sufficiently close to these lower bounds.
In the limit \(\theta\downarrow 3/p\) and \(\beta\downarrow 1/p\), the last inequality becomes
\begin{equation*}
\frac32 q>\frac32+\frac{5}{2p},
\end{equation*}
which is equivalent to the left inequality of \eqref{eq:sw-cond}.
Thus suitable \(\theta\) and \(\beta\) exist.

The remaining condition in \eqref{eq:sw-cond},
\begin{equation*}
0\le \frac{4p-1}{pq-1}-3,
\end{equation*}
is not used in the above choice of \(\theta\) and \(\beta\).
It only indicates that the present case lies outside the strong--strong region covered by Corollary \ref{cor:strong-strong}.
\end{proof}

\begin{remark}
The condition \(\beta p>1\) is essential in the above contraction argument.
It is the condition that allows the cone weight required in the damped wave estimate to be supplied by the decay of \(|v|^p\).
Therefore, in Corollaries \ref{cor:strong-strong} and \ref{cor:weak-strong}, the statement should be read as the existence of a suitable \(\beta\in(1/p,1)\), not as a result for every \(\beta\in(0,1)\).
\end{remark}

\noindent\textbf{Acknowledgements}\quad V. G. is partially supported by the Gruppo Nazionale per l'Analisi Matematica, the Institute of Mathematics and Informatics of the Bulgarian Academy of Sciences, and the Top Global University Project at Waseda University.
K. K. is partially supported by JSPS KAKENHI Grant Numbers JP24K16955 and JP24KJ0010.


\begin{thebibliography}{99}

\bibitem{Agemi-Kurokawa-Takamura}
R. Agemi, Y. Kurokawa and H. Takamura, \textit{Critical curve for \(p\)-\(q\) systems of nonlinear wave equations in three space dimensions}, J. Differential Equations, \textbf{167} (2000), no. 1, 87--133.

\bibitem{Chen-Liu}
W. Chen and Y.Liu, \textit{Large time behavior for the hyperbolic--parabolic coupled system with the regularity-loss structure}, Journal of Evolution Equations, \textbf{25} (2025), 61 pages.

\bibitem{Chen-Reissig}
W. Chen and M. Reissig, \textit{Blow-up of solutions to Nakao's problem via an iteration argument}, J. Differential Equations, \textbf{275} (2021), 733--756.

\bibitem{Chill-Haraux}
R. Chill and A. Haraux, \textit{An optimal estimate for the difference of solutions of two abstract evolution equations}, J. Differential Equations, \textbf{193} (2003), no. 2, 385--395.

\bibitem{DelSanto-Georgiev-Mitidieri}
D. Del Santo, V. Georgiev and E. Mitidieri, \textit{Global existence of the solutions and formation of singularities for a class of hyperbolic systems}, In: \textit{Geometrical Optics and Related Topics}, Progr. Nonlinear Differential Equations Appl., \textbf{32}, Birkh\"auser Boston, Boston, MA, 1997, 117--140.

\bibitem{Escobedo-Herrero}
M. Escobedo and M. A. Herrero, \textit{Boundedness and blow up for a semilinear reaction-diffusion system}, J. Differential Equations, \textbf{89} (1991), no. 1, 176--202.

\bibitem{Fujita}
H. Fujita, \textit{On the blowing up of solutions of the Cauchy problem for \(u_t=\Delta u+u^{1+\alpha}\)}, J. Fac. Sci. Univ. Tokyo Sect. I, \textbf{13} (1966), 109--124.

\bibitem{Georgiev-Heiming-Kubo}
V. Georgiev, Ch. Heiming and H. Kubo, \textit{Supercritical semilinear wave equation with non-negative potential}, Comm. Partial Differential Equations, \textbf{26} (2001), no. 11--12, 2267--2303.

\bibitem{Georgiev-Kita}
V. Georgiev and K. Kita, \textit{Weighted \(L^\infty\)-estimates for solutions of the damped wave equation in three space dimensions and its application}, Nonlinear Analysis, \textbf{260} (2025), Paper No. 113850.

\bibitem{Georgiev-Takamura-Zhou}
V. Georgiev, H. Takamura and Y. Zhou, \textit{The lifespan of solutions to nonlinear systems of a high-dimensional wave equation}, Nonlinear Anal., \textbf{64} (2006), no. 10, 2215--2250.

\bibitem{Ikeda-Ogawa}
M. Ikeda and T. Ogawa, \textit{Lifespan of solutions to the damped wave equation with a critical nonlinearity}, J. Differential Equations, \textbf{261} (2016), no. 3, 1880--1903.

\bibitem{Ikeda-Sobajima-Wakasa}
M. Ikeda, M. Sobajima and K. Wakasa, \textit{Blow-up phenomena of semilinear wave equations and their weakly coupled systems}, J. Differential Equations, \textbf{267} (2019), no. 9, 5165--5201.

\bibitem{Ikeda-Wakasugi}
M. Ikeda and Y. Wakasugi, \textit{A note on the lifespan of solutions to the semilinear damped wave equation}, Proc. Amer. Math. Soc., \textbf{143} (2015), no. 1, 163--171.

\bibitem{Ikehata-Nishihara}
R. Ikehata and K. Nishihara, \textit{Diffusion phenomenon for second order linear evolution equations}, Studia Math., \textbf{158} (2003), no. 2, 153--161.

\bibitem{Ikehata-Tanizawa}
R. Ikehata and K. Tanizawa, \textit{Global existence of solutions for semilinear damped wave equations in \(\mathbb R^N\) with noncompactly supported initial data}, Nonlinear Anal., \textbf{61} (2005), no. 7, 1189--1208.

\bibitem{John}
F. John, \textit{Blow-up of solutions of nonlinear wave equations in three space dimensions}, Manuscripta Math., \textbf{28} (1979), no. 1--3, 235--268.

\bibitem{Kita-Kusaba}
K. Kita and R. Kusaba, \textit{A remark on the blowing up of solutions to Nakao's problem}, J. Math. Anal. Appl., \textbf{513} (2022), no. 1, Paper No. 126199.

\bibitem{Kurokawa-Takamura-Wakasa}
Y. Kurokawa, H. Takamura and K. Wakasa, \textit{The blow-up and lifespan of solutions to systems of semilinear wave equation with critical exponents in high dimensions}, Differential Integral Equations, \textbf{25} (2012), no. 3--4, 363--382.

\bibitem{Hsiao-Liu}
L. Hsiao and  T.-P. Liu, \textit{Convergence to nonlinear diffusion waves for solutions of a system of hyperbolic conservation laws with damping}, Communications in Mathematical Physics, \textbf{143} (1992), no. 3, 599--605.

\bibitem{Kawashima}
S. Kawashima, \textit{Large-time behaviour of solutions to hyperbolic-parabolic systems of conservation laws and applications}, Proceedings of the Royal Society of Edinburgh, Section A: Mathematics, \textbf{106} (1987), 169-194.

\bibitem{Li-Zhou}
T.-T. Li and Y. Zhou, \textit{Breakdown of solutions to \(\Box u+u_t=|u|^{1+\alpha}\)}, Discrete Contin. Dyn. Syst., \textbf{1} (1995), no. 4, 503--520.

\bibitem{Lindblad}
H. Lindblad, \textit{Blow-up for solutions of \(\Box u=|u|^p\) with small initial data}, Comm. Partial Differential Equations, \textbf{15} (1990), no. 6, 757--821.

\bibitem{Matsumura}
A. Matsumura, \textit{On the asymptotic behavior of solutions of semi-linear wave equations}, Publ. Res. Inst. Math. Sci., \textbf{12} (1976/77), no. 1, 169--189.

\bibitem{Narazaki}
T. Narazaki, \textit{Global solutions to the Cauchy problem for the weakly coupled system of damped wave equations}, Discrete Contin. Dyn. Syst. Suppl. (2009), 592--601.

\bibitem{Narazaki2004}
T. Narazaki, \textit{\(L^p\)-\(L^q\) estimates for damped wave equations and their applications to semi-linear problem}, J. Math. Soc. Japan, \textbf{56} (2004), no. 2, 585--626.

\bibitem{Narazaki-Nishihara}
T. Narazaki and K. Nishihara, \textit{Asymptotic behavior of solutions for the damped wave equation with slowly decaying data}, J. Math. Anal. Appl., \textbf{338} (2008), no. 2, 803--819.

\bibitem{Nishihara1997}
K. Nishihara, \textit{Asymptotic behavior of solutions of quasilinear hyperbolic equations with linear damping}, J. Differential Equations, \textbf{137} (1997), no. 2, 384--395.

\bibitem{Nishihara2003}
K. Nishihara, \textit{\(L^p\)-\(L^q\) estimates of solutions to the damped wave equation in 3-dimensional space and their application}, Math. Z., \textbf{244} (2003), no. 3, 631--649.

\bibitem{Nishihara-Wakasugi}
K. Nishihara and Y. Wakasugi, \textit{Critical exponent for the Cauchy problem to the weakly coupled damped wave system}, Nonlinear Anal., \textbf{108} (2014), 249--259.

\bibitem{Shizuta-Kawashima}
Y. Shizuta and S. Kawashima, \textit{Systems of equations of hyperbolic-parabolic type with applications to the discrete Boltzmann equation}, Hokkaido Mathematical Journal, \textbf{14}, (1985), no. 2,  249 -- 275.

\bibitem{Strauss}
W. A. Strauss, \textit{Nonlinear scattering theory at low energy}, J. Funct. Anal., \textbf{41} (1981), no. 1, 110--133.

\bibitem{Sun-Wang}
F. Sun and M. Wang, \textit{Existence and nonexistence of global solutions for a nonlinear hyperbolic system with damping}, Nonlinear Anal., \textbf{66} (2007), no. 12, 2889--2910.

\bibitem{Todorova-Yordanov}
G. Todorova and B. Yordanov, \textit{Critical exponent for a nonlinear wave equation with damping}, J. Differential Equations, \textbf{174} (2001), no. 2, 464--489.

\bibitem{Wakasugi}
Y. Wakasugi, \textit{A note on the blow-up of solutions to Nakao's problem}, In: \textit{New Trends in Analysis and Interdisciplinary Applications}, Birkh\"auser/Springer, 2017, 545--551.

\bibitem{Wirth}
J. Wirth, \textit{Diffusion phenomena for partially dissipative hyperbolic systems}, Journal of Mathematical Analysis and Applications, \textbf{414}, (2014), 666--677.

\bibitem{Yang-Milani}
H. Yang and A. Milani, \textit{On the diffusion phenomenon of quasilinear hyperbolic waves}, Bull. Sci. Math., \textbf{124} (2000), no. 5, 415--433.

\bibitem{Zhang}
Q. S. Zhang, \textit{A blow-up result for a nonlinear wave equation with damping: the critical case}, C. R. Acad. Sci. Paris S\'er. I Math., \textbf{333} (2001), no. 2, 109--114.

\end{thebibliography}
\end{document}